# Structures in P based on Properties of Semigroup and Arithmetical Sequence H = (±3·2; 1)

**Author: Michael H. Hebert**

**Abstract**
This paper presents results on structures in P based on tools developed from subjects of elementary number theory. Key findings are: The arithmetical sequence H = (±3·2; 1) is in Z the smallest superset of P \ {3, 2}. H is a semigroup. A revised definition of P. Unique Gestalt of p in Z. The prime number lattice packing $H^n$. The geometrical locus in $H^2$ of the family of solutions of: - the set of prime twins, - the set of PRACHAR prime twins, - in $H^2$, $H^3$ family of solutions of the GOLDBACH conjunction. Partition of H in $p^2$-intervals. Prime numbers in $p^2$-intervals. Infinity of the set of prime twins. Verification of the GOLDBACH conjunction.

**Zusammenfassung**
Diese Arbeit zeigt Ergebnisse zu Strukturen in P auf Basis von Werkzeugen, die aus Begriffen der elementaren Zahlentheorie entwickelt wurden. Wesentliche Ergebnisse sind: Die arithmetische Folge H = (±3·2; 1) ist in Z die kleinste Obermenge von P \ {3, 2}. H ist eine Halbgruppe. Eine überarbrbeitete Defintion von P. Einheitliche Gestalt von p in Z. Das Primzahl-Gitter $H^n$. Die geometrischen Örter in $H^2$ der Lösungsmenge: - der Primzahlzwillinge, - der PRACHAR Primzahlzwillinge, - in $H^2$, $H^3$ der GOLDBACH-Vermutung. Intervall-Einteilung von H in $p^2$-Intervalle. Primzahlen in $p^2$-Intervallen. Unendlichkeit der Menge der Primzahlzwillinge. Bestätigung der GOLDBACH-Vermutung.

**Dank** gilt Prof. JAN-CHRISTOPH SCHLAGE-PUCHTA, Universität Rostock, für viele hilfreiche Hinweise.

**Inhalt**

**Klassische Sätze, die hier zur Anwendung kommen**
**Satz I.** Fundamentalsatz der elementaren Zahlentheorie: „Jede natürliche Zahl a > 1 ist als Produkt von Primzahlen darstellbar: a = $p_1\, p_2 \dots p_m$. Die Darstellung ist, abgesehen von der Reihenfolge der Faktoren, eindeutig.“ ([2], Seite 10.)
**Satz II.** DIRICHLET: Schubfächer-Prinzip: „Verteilt man n Dinge auf m Klassen und ist m < n, so kommen in irgendeiner Klasse mindestens zwei Dinge vor.“ ([2], Seite 6.)
**Satz III.** DIRICHLET (1837): Primzahlen in arithmetischen Folgen: „Jede arithmetische Folge {a + k·b I a ∈ Z, b ∈ N, k = 0, 1, 2, …}, in der das Anfangsglied a und die Differenz b zueinander teilerfremd sind, enthält unendlich viele Primzahlen.“ ([3], Seite 140, sinngemäß.) oder ([4], Seite 27.)
**Satz IV.** HADAMARD (1896), DE LA VALLEÉ POUSSIN (1896): Primzahlsatz: „Mit $x \to \infty$ ist $\lim \pi(x) \cdot (x / \log x)^{-1} = 1$.“ ([3], Seite 12.)

# Strukturen in P

## § 1. Die kleinste Obermenge von P

Arithmetische Zahlenfolgen

Gegeben sei die mit N bezeichnete *Folge* der *natürlichen Zahlen* {0, 1, 2, …} := N. Dann ist auch die mit Z bezeichnete Folge der *ganzen Zahlen* {…, -2, -1, 0, 1, 2, …} := Z gegeben. Der arithmetische und algebraische Umgang mit N und Z wird vorasgesetzt. - Eine *arithmetische Zahlenfolge*, aZ, ist eine Reihung unendlich vieler natürlicher Zahlen {a, a+d, a+2d, …} := (+d; a) mit dem *Anfangsglied* a und der *konstanten Differenz* d = $a_{n+1} – a_n$. Um zum nächsten Glied der Folge zu kommen, wird bei a beginnend d addiert. Beispiel: N mit a = 0, d = 1. Durch a und d ist die aZ bestimmt. - Wir ergänzen die Definition um ‚-' und ‚±' vor d. Bei -d wird, um zum nächsten Glied zu kommen, subtrahiert. Bei ±d wird addiert und subtrahiert; die Folge ragt dann in den negativen und in den positiven Bereich von Z. Beispiel: Z mit a = 0. Man sagt: bei ‚+d' hat die Folge 1 *Ast*, der nach +∞ strebt. Bei ‚-d' hat sie 1 Ast, der nach -∞ strebt. Bei ±d hat sie 2 Äste und strebt von -∞ nach +∞. Das Streben der aZ geschieht ohne *Häufungspunkte*.

Herkömmliche Definition der Primzahlen

Sie lautet:

**Definition 1.1.** Primzahl (herkömmliche Definition):

*Primzahl p: ≠1, positiv, ganz, nur durch 1 und durch sich selbst teilbar.*

Wir nennen diese Definition ‚herkömmlich'; denn wir werden eine einfachere vorschlagen. Die *Folge der Primzahlen* {2, 3, 5, …} wird mit P bezeichnet; einzelne Primzahlen mit p oder kleinen lateinischen Buchstaben. Bereits EUKLID bewies:

**Satz 1.1.** Unendlichkeit P:

*Es gibt unendlich viele Primzahlen* (EUKLID).

Zur herkömmlichen Definition der Primzahlen fällt z.B. auf: Die Zahl 1 wird von P ausgeschlossen und erhält besondere Behandlung. Unsere Definition siehe Def. 1.3.

Kleinste Obermenge von P

In der Def. 1.1. ist N die *Obermenge* von P. Das führt zu der Frage, ob es für P andere, kleinere Obermengen gibt - und welche die kleinste ist. Wir suchen die kleinste Obermenge von P, lassen uns aber von der Def. 1.1. nicht einengen. Wir folgen nur dem mathemat. Prinzip, stets mit den geringsten Mitteln auszukommen. Die Suche geschieht in zwölf Schritten:

(1) X = {…, $x_n$, …} sei ein Modell der gesuchten Menge.

(2) Die Glieder $x_n$ von X seien ganze Zahlen: $x_n \in$ Z. Wir verlagern also unsere *Grundmenge* von N nach Z.

(3) Z ist eine aZ und ist eine Obermenge von P. Wir lassen zu, dass p < 0 sein darf.

(4) Z ist für unsere Zwecke ‚zu groß'.

(5) Wir ‚verkleinern' Z, etwa durch Streichungen.

(6) Die Konstruktion von X als aZ soll erhalten bleiben.

(7) Die Frage nach anderen Obermengen von P verschärfen wir zur Frage nach kleineren aZ = (±d; a) = X.

(8) Von den in Z definierten *Verknüpfungen Mult*., *Add*., *Subtr*. u. *Div*. können wir nur die Multiplikation beibehalten. Denn die anderen Verknüpfungen können die aZ verlassen.

(9) Wir suchen also eine aZ mit der Gestalt X = (±d; a) und mit der Multiplikation als zweistelliger Verknüpfung.

(10) Sei in (±d; a) d *gerade* und a *ungerade*. Dann und nur dann sind, wie für Primzahlen erforderlich, alle Glieder ungerade. Seien x u. y beliebige Glieder von (±d; a), etwa x = d·m+a und y = d·n+a. Dann ist x·y = d·(...) + $a^2$. Also x·y ∈ (±d; a) dann und nur dann, wenn a = $a^2$, d.h. a = 1. D.h.: *Ist eine aZ = (±d; a) multiplikativ, so ist a = 1.*

(11) Die Multiplikation hat in (±d; 1) die Eigenschaften: 1) Eindeutigkeit: Je zwei Elemente a, b ∈ (±d; 1) besitzen eindeutig ein Element c ∈ (±d; 1) mit c = a·b; wegen der Eindeutigkeit der Multiplikation und der Gestalt der Faktoren. 2) Assoziativität: Für je drei Elemente a, b, c ∈ (±d; 1) gilt: (a·b)·c = a·(b·c); wegen der Assoziat. der Multiplikation. 3) Eins-Element: a·e = e·a = a, für jedes a ∈ (±d; 1) und mit der Zahl e = 1 als Eins-Element. 4) Kommutativität: a·b = b·a, für jede a, b ∈ (±d; 1); wegen der Kommutativität der Multiplikation. 5) Inverses Element: Ein inverses Element existiert nicht generell, weil die Division nur eingeschränkt gilt. Also: 1), 2): (d; 1) ist Halbgruppe. 3): (d; 1) hat Eins-Element. 4): (d; 1) ist kommutativ. 5): (d; 1) ist nicht Gruppe: D.h.: *Ist eine aZ = (±d; 1), so ist sie eine kommutative Halbgruppe mit Eins-Element.*

(12) X soll *fast alle* p enthalten, d.h. alle mit endlich vielen Ausnahmen. Sie soll die ‚kleinste' Obermenge von P sein – (Bemerkung: Angesichts dieser Forderungen ist N nicht die kleinste Obermenge von P.) Wenn die gesuchte Obermenge X = (±d; a) eine aZ sein soll, dann gibt es nur die Parameter a und d zur Unterscheidung der Kandidaten. Bereits festgestellt ist a = 1. Klar ist: Je größer d, desto kleiner die Folge. Eine systematische Suche nach der kleinsten Obermenge geschieht also über wachsendes d beginnend bei d = 2. Es sind aZ-typische Rechenregeln zu beachten: 1.) (a,d)=1 *teilerfremd*, 2.) a < d, 3.) (d; a) = (d; d-a) . Notwendig: d > 1; denn dann kommt jede Primzahl in X nur einmal, d.h. genau einmal vor. Also gibt es nur die folgenden Kombinationen von d und a, die sinnvolle Lösungen erwarten lassen. Wir zeigen alle Kombinationen. Dabei wächst d von 2 nach 6:
X = (±2; 1) = {…, -5, -3, -1, 1, 3, 5, …}, enthält alle p außer 2, nicht klein genug,
oder X = (±3; 1) = -(±3; 2) = {…, -8, -5, -2, 1, 4, 7, …}, enthält alle p außer 3, nicht klein genug,
oder X = (±4; 1) = -(±4; 3) = {…, -11, -7, -3, 1, 5, 9,…}, enthält alle p außer 2, nicht klein genug,
oder X = (±5; 1) = -(±5; 4) = {… , -14, -9, -4, 1, 6, 11,…}, es fehlen: 2, 3, 5, …
oder X = (±6; 1) = -(±6; 5) = {…, -17,-11, -5, 1, 7, 13, 19, …}, enthält alle p außer 2 und 3.
Ab d = 7 entstehen größere Lücken bei den p. Die Vorgehensweise beweist außerdem, dass eine ‚Minimal'-Lösung gefunden wurde. Die gesuchte ‚kleinste' Obermenge von P \ {3, 2} ist X := (±3·2; 1); sie wird H genannt:
**Definition 1.2.** Die Folge H:
*H := (±3·2; 1)* = {…, -17,-11, -5, 1, 7, 13, 19, …}
**Satz 1.2.** H kleinste Obermenge:
*H ist in Z die kleinste aZ, die Obermenge von P \ {3, 2} ist.*
**Satz 1.3.** H Halbgruppe:
*H ist eine kommutative Halbgruppe mit 1-Element.*
H ist d i e Ausgangs-Folge der vorliegenden Arbeit.

Offenbar gilt auch:
**Satz 1.4.** Gestalt der Glieder in H:
*In H haben a l l e Glieder, also auch die Primzahlen, die Gestalt 3·2·k +1 mit k ∈ Z.*

Unsere vereinfachte Definition der Primzahlen kann nun lauten:
**Definition 1.3.** Primzahl:
*Primzahl p: Glied von H, nur durch 1 und durch sich selbst teilbar. Außerdem: 2 ∈P, 3 ∈P.*
Bemerkung: Die 1 ist n i c h t ausgeschlossen, es gibt für sie keine besondere Behandlung und die Primzahlen liegen auch im n e g a t i v e n Ast von Z.

Primzahlen in den Ästen von H
Wir teilen die Folge H in ihre beiden Äste und klappen den negativen Ast ins Positive durch Multiplikation mit -1. Das führt zur Def. der aZ $H_+$ und $H_-$:
**Definition 1.4.** Äste von H:
*$H_+$ := (+3·2; +1) = {1, 7, 13, …},*
*$H_-$ := (+3·2; -1) = {-1, 5, 11, …}.*

Der folgende Satz 1.5. gilt insbes. auch für Primzahlen:
**Satz 1.5.** *Gestalt in $H_+$ und $H_-$:*
*Für x ∈ $H_+$ ex. k ∈ N mit x = 3·2·k +1.*
*Für x ∈ $H_-$ ex. k ∈ N mit x = 3·2·k -1.*

Die Folgen $H_+$ und $H_-$ sind bereits bei DIRICHLET anzutreffen:
**Satz 1.6.** Satz von DIRICHLET über Primzahlen in der arithmetischen Progression (Satz III):
*$H_+$ enthält ∞ viele Primzahlen, und zwar die mit der Gestalt p = 3·2 ·k + 1, k ∈ N.*
*$H_-$ enthält ∞ viele Primzahlen, und zwar die mit der Gestalt p = 3·2 ·k - 1, k ∈ N.*

## § 2. Das Primzahl-Gitter.
In $H_-$ und $H_+$ ist die Gestalt der Glieder x = 3·2·s±1. s heiße *Träger*. $H_-$ und $H_+$ haben dieselben Träger: die natürlichen Zahlen. $S_-$ und $S_+$ seien die Träger-Folgen von $H_-$ und $H_+$; es ist $S_-$ = $S_+$ = N. Der Träger zeigt die Position eines Gliedes in $H_-$ bzw. $H_+$ in den beiden ‚symmetrischen' Folgen. Die Positionen der Primzahlen und Vielfachen in $H_-$ und $H_+$ werden bewahrt; das ist f u n d a m e n t a l für unseren Ansatz. – Mit der Folge H bzw. mit ihren Ästen spannen wir das Primzahl-Gitter auf:

**Definition 2.1.** *Primzahl-Gitter 2-dimensional* $:= H \times H = H^2 = \{(x, y) \mid x, y \in H\}$, *…3-dimensional:* $H \times H \times H: = H^3 = \{(x, y, z) \mid x, y, z \in H\}$, *usw.*

Tabelle 2.1. zeigt Anfangswerte von $H^2$. In $H^2$ bewahren die Primzahlen und die Vielfachen die Positionen aus H. Wir verstehen o.B.d.A. $H^2$ als *rechtwinkliges kartesisches Koordinaten-System*. *Ursprung* (+1,+1). $H^2$ hat die *Quadranten*: 1. $H_+ \times H_+$, 2. $H_+ \times H_-$, 3. $H_- \times H_-$ und 4. $H_- \times H_+$. Die Positionierung der Primzahlen und *Vielfachen* generiert in den Quadranten eine Rechteck- und Streifen-Struktur. pp: Punkt hat zwei Primzahl-Komponenten. mm: Punkt hat zwei Vielfachen-Komponenten. pm oder mp, in der Tabelle · : Punkt hat eine Primzahl- und eine Vielfachen-Komponente. Andere Kombinationen gibt es nicht. Im 1. und 3. Quadranten besteht Symmetrie zur Winkelhalbierenden der Quadranten. Nicht so im 2. und 4. Quadranten. Man könnte auch (-1, -1) als Ursprung wählen. Aber wir verfahren nicht so, weil die Primzahl-Quadrate, die für die spätere Definition der Relevanz-Intervalle gebraucht werden, Glieder von $H_+$ sind. Es gilt:

**Satz 2.1.** *Eigenschaften des Primzahl-Gitters $H^2$: 1) Sein Koordinaten-Ursprung ist der Punkt (+1, +1). 2) Es wird aufgespannt von den Ästen von H. 3) Seine Träger-Folgen $S_-$ und $S_+$ sind einander gleich und sind = N. 4) Es hat die vier Quadranten 1. $H_+ \times H_+$, 2. $H_+ \times H_-$, 3. $H_- \times H_-$ und 4. $H_- \times H_+$. 5) Symmetrie besteht nur zwischen dem 1. und 3. Quadranten, d. h. zur Winkelhalbierenden der Quadranten. 6) Es besitzt eine bestimmte geometrische Struktur: Primzahlen und Vielfache generieren Rechtecke und Streifen. Alles analog in $H^3, \ldots, H^n$.*

Bemerkung: Die Struktur, die das Primzahl-Gitter erzeugt, wird nur sichtbar, wenn man Z als Obermenge von P wählt.

2. Quadrant — H+ ↑ — 1. Quadrant

| | | | | | | | | | | | | | | | | | | | | | | |
|---|---|---|---|---|---|---|---|---|---|---|---|---|---|---|---|---|---|---|---|---|---|---|
| | pp | pp | pp | pp | . | pp | pp | pp | pp | pp | 61 | pp | pp | pp | . | pp | pp | pp | . | . | pp | |
| | | | | . | mm | . | . | . | . | . | 55 | . | . | . | mm | . | . | . | mm | | mm | |
| | . | . | . | . | mm | . | . | . | . | . | 49 | . | . | . | mm | . | . | . | mm | | mm | |
| | pp | pp | pp | pp | . | pp | pp | pp | pp | pp | 43 | pp | pp | pp | . | pp | pp | pp | . | . | pp | |
| | pp | pp | pp | pp | . | pp | pp | pp | pp | pp | 37 | pp | pp | pp | . | pp | pp | pp | . | . | pp | |
| | pp | pp | pp | pp | . | pp | pp | pp | pp | pp | 31 | pp | pp | pp | . | pp | pp | pp | . | . | pp | |
| | . | . | . | . | mm | . | . | . | . | . | 25 | . | . | . | mm | . | . | . | mm | | mm | |
| | pp | pp | pp | pp | . | pp | pp | pp | pp | pp | 19 | pp | pp | pp | . | pp | pp | pp | . | . | pp | |
| | pp | pp | pp | pp | . | pp | pp | pp | pp | pp | 13 | pp | pp | pp | . | pp | pp | pp | . | . | pp | |
| | pp | pp | pp | pp | . | pp | pp | pp | pp | pp | 7 | pp | pp | pp | . | pp | pp | pp | . | . | pp | |
| $H_-$ ← | -59 | -53 | -47 | -41 | -35 | -29 | -23 | -17 | -11 | -5 | 1 | 7 | 13 | 19 | 25 | 31 | 37 | 43 | 49 | 55 | 61 | → $H_+$ |
| | pp | pp | pp | pp | . | pp | pp | pp | pp | pp | -5 | pp | pp | pp | . | pp | pp | pp | . | . | pp | |
| | pp | pp | pp | pp | . | pp | pp | pp | pp | pp | -11 | pp | pp | pp | . | pp | pp | pp | . | . | pp | |
| | pp | pp | pp | pp | . | pp | pp | pp | pp | pp | -17 | pp | pp | pp | . | pp | pp | pp | . | . | pp | |
| | pp | pp | pp | pp | . | pp | pp | pp | pp | pp | -23 | pp | pp | pp | . | pp | pp | pp | . | . | pp | |
| | pp | pp | pp | pp | . | pp | pp | pp | pp | pp | -29 | pp | pp | pp | . | pp | pp | pp | . | . | pp | |
| | | | | | m . | . | . | . | | . | -35 | | . | | mm | . . | . . | | | mm | mm | |
| | pp | pp | pp | pp | . | pp | pp | pp | pp | pp | -41 | pp | pp | pp | . | pp | pp | pp | . | . | pp | |
| | pp | pp | pp | pp | . | pp | pp | pp | pp | pp | -47 | pp | pp | pp | . | pp | pp | pp | . | . | pp | |
| | pp | pp | pp | pp | . | pp | pp | pp | pp | pp | -53 | pp | pp | pp | . | pp | pp | pp | . | . | pp | |
| | pp | pp | pp | pp | . | pp | pp | pp | pp | pp | -59 | pp | pp | pp | . | pp | pp | pp | . | . | pp | |

3. Quadrant — ↓ $H_-$ — 4. Quadrant

Primzahl-Gitter $H \times H = (\pm 3 \cdot 2 \cdot s; +1) \times (\pm 3 \cdot 2 \cdot s; +1)$, Anfangswerte
Tabelle 2.1.

## § 3. Geometrischer Ort der Lösungsmenge der 2-gliedrigen GOLDBACH-Vermutung.

<u>Klassische Vermutungen</u>

Wir ermitteln die geometrischen Örter im Primzahl-Gitter der Lösungsmengen der 2-dimensionalen GOLDBACHschen Vermutung (2GH), der Primzahlzwillinge (PT), der PRACHAR Primzahlzwillinge (PPT) und der (3-dimensionalen) GOLDBACHschen Vermutung (3GH). Die Definitionen dieser Probleme sind z.B. bei PRACHAR zu finden. Es seien $p, p_i, p_{i-1}, q \in P$ und $i \in N$.

(2GH): „GOLDBACH vermutete 1742, dass sich jede gerade Zahl > 2 als Summe zweier Primzahlen darstellen lasse.“ [1, S. 51].

(PT): „Seit langer Zeit wird vermutet, dass es unendlich viele i gibt, für die $d_i = p_i - p_{i-1} = 2$ ist. “ [1, S.155].

(PPT): „Es ist bis heute nicht bekannt, ob es $\infty$ viele p gibt, für die $\frac{1}{2} \cdot (p - 1)$ wieder prim ist.“ [1, S. 51, Fußnote].

(3GH): „Jede ungerade Zahl > 7 ist Summe dreier Primzahlen." [3, S. 172].

Zu beweisen

Zu beweisende Vermutungen:

(2GH): "Zu i, $1 < i \in N$, existieren p, $q \in P$, so dass $p + q = 2 \cdot i$."

(PT): "Die Menge $\{(p, q)|\ p, q \in P, p = q + 2\}$ ist unendlich."

(PPT): "Die Menge $\{(p, q)|\ p, q \in P, p = ½ \cdot (q - 1)\}$ ist unendlich." Wir nennen (PPT) 'PRACHAR Primzahlzwillinge'.

Harmonisierung

Eine Harmonisierung dieser klassischen linearen Formulierungen ist angebracht: dazu schreiben wir diese p-q-Gleichungen als Geraden-Gleichungen der x-y-Ebene:

(2GH): $x = -1 \cdot y + 2 \cdot i$, (PT): $x = 1 \cdot y + 2$, (PPT): $x = 2 \cdot y + 1$.

Die Gleichungen motivieren, die Lösungen als Punkte von Geraden einer Ebene zu verstehen. Wir suchen eine geeignete Ebene. ‚Lösung' (p, q) - bzw. (x, y) - heißt ein geordnetes Paar (p, q) von Primzahlen, das (2GH) oder (PT) oder (PPT) erfüllt. p u. q heißen Lösungskomponenten.

(2GH)

Zu Beginn ist nicht bekannt, welche Primzahlen Lösungskomponenten sind. Daher ist N die Ausgangsmenge für (2GH). Lösungen mit den Lösungskomponenten 2 oder 3 nennen wir trivial und schließen sie von der Lösungsmenge aus.

$\{\ldots, 2 \cdot n, \ldots\}$, $n \in N$, seien die geraden Zahlen. Fälle $2 \cdot n=0$ und $2 \cdot n=2$ sind nicht Gegenstand von (2GH). Fälle 4=2+2, 6=3+3, 8=3+5 haben nur triviale Lösungen. Wir betrachten die geraden Zahlen e > 8. Sie bilden die aZ $e_n$ := (+2; +10) = {10, 12, 14, …}. In $e_n$ kommt jede gerade Zahl e > 8 genau ein Mal vor. Wir beziehen die Glieder auf ihre Reste (mod $3 \cdot 2$). Es gilt trivial:

**Satz 3.1.** *__Die Reste -2, 0, +2 (mod 3·2) der geraden Zahlen__: $x \in (2; 10) \Rightarrow$ s existiert mit $x := 3 \cdot 2 \cdot s -2$ oder $x := 3 \cdot 2 \cdot s \pm 0$ oder $x = :3 \cdot 2 \cdot s + 2$, $1 < s \in N$. D.h.: $x \equiv -2$ oder $\equiv 0$ oder $\equiv +2$ (mod $3 \cdot 2$).*

**Satz 3.2.** *__Zerlegung von (+2; +10) in die Folgen $E_{-2}$, $E_0$ und $E_{+2}$:__ $E_{+2}$ := (3·2; 14) = {14, 20, 26, …} $\equiv +2$ (mod 3·2), $E_0$ := (3·2; 12) = {12, 18, 24, …} $\equiv 0$ (mod 3·2), $E_{-2}$ := (3·2; 10) = 10, 16, 22, … $\equiv -2$ (mod 3·2). $E_{-2}$, $E_0$ und $E_{+2}$ sind zueinander fremd und vereinigt gleich $e_n$. In $e_n$ kommt jede gerade Zahl e > 8 genau 1 Mal vor.*

Wir betrachten auch die Folgen $S(E_{+2})$, $S(E_{-2})$ und $S(E_0)$ der Träger von $E_{+2}$, $E_{-2}$ und $E_0$. Sie sind einander gleich und sind gleich N \ {0, 1}. Sie gewährleisten die Vollständigkeit der Lösungsmenge:

**Satz 3.3.** __Die Träger-Folgen:__ *$S(E_{+2}) = S(E_0) = S(E_{-2}) = N \setminus \{0, 1\}$.*

Wir betrachten x + y mit x, y > 4 und x, $y \in |H|$, also der Gestalt $3 \cdot 2 \cdot s \pm 1$, und ihre Kongruenzen (mod 3·2). x + y ist gerade, denn x und y sind ungerade. Es gibt genau drei Fälle:

**Satz 3.4.** __Fälle__: *Fall 1: $x, y \in H_+ \Rightarrow x+y \in E_{+2}$. Fall 2: $x, y \in H_- \Rightarrow x+y \in E_{-2}$. Fall 3: $x \in H_+$, $y \in H_-$, oder umgekehrt $\Rightarrow x+y \in E_0$.*

Fall 1: $x+y = 3 \cdot 2 \cdot a+1 + 3 \cdot 2 \cdot b+1 = 3 \cdot 2 \cdot (a+b)+2 \equiv +2$ (mod 3·2). Fall 2: $x+y = 3 \cdot 2 \cdot a-1 + 3 \cdot 2 \cdot b-1 = 3 \cdot 2 \cdot (a+b)-2 \equiv - 2$ (mod 3·2). Fall 3: $x+y = 3 \cdot 2 \cdot a+1 + 3 \cdot 2 \cdot y-1 = 3 \cdot 2 \cdot (a+b) \equiv 0$ (mod 3·2). Mit der Beschränkung auf Primzahlen dürfen wir sagen:

**Satz 3.5.** *__Die Lösungen von (2GH)__: Wenn p + q = e, $(3 < p, q) \wedge (p, q \in P) \wedge (e \in E_{+2}$ oder $E_{-2}$ oder $E_0)$, so (p, q) Lösung von (2GH).*

Wir teilen Satz 3.5. nach Kongruenzen auf und zeigen die Träger:

**Satz 3.6.** *__Gestalt der Lösungen__: Fall 1: $3 \cdot 2 \cdot s_p+1 + 3 \cdot 2 \cdot s_q+1 = 3 \cdot 2 \cdot s_e \equiv +2$ (mod 3·2), Fall 2: $3 \cdot 2 \cdot s_p -1 + 3 \cdot 2 \cdot s_q -1 = 3 \cdot 2 \cdot s_e \equiv -2$ (mod 3·2), Fall 3: $3 \cdot 2 \cdot s_p \pm 1 + 3 \cdot 2 \cdot s_q \pm 1 = 3 \cdot 2 \cdot s_e \equiv 0$ (mod 3·2).*

$s_p$ und $s_q$ sind die Träger von p und q und $s_e$ ist der von e. Wir dürfen für e die Träger-Darstellung verwenden. Sie ist nicht reserviert für die Gestalt $3 \cdot 2 \cdot r \pm 1$; sie hängt nur ab von der Relation zu 3·2. Nach Satz 3.1. und Satz 3.6. folgt:

**Satz 3.7.** __Träger-Schreibweise__: *Einheitliche Darstellung der Form der Lösungen mittels der Träger-Schreibweise: $s_p + s_q = s_e$.*

Satz 3.1. bis Satz 3.7.sind evident.

Die Lösungsmenge im Allgemeinen

Die Lösungen werden notiert in H×H. Für die (2GH) sind Primzahlen positiv. Deshalb verwenden wir das Gitter in der Form $|H| \times |H|$. Die Gitter-Punkte seien besetzt von den Summen ihrer Koordinaten. Bei den Gitter-Punkten gibt es die Fälle pp, pm, mp und mm. pp: beide Komponenten Primzahlen und es handelt sich um eine Lösung. Der Anfang der Lösungsmenge erscheint explizit in Tabelle 3.1. Zur Verdeutlichung sind nur die pp-Fälle notiert.

Tabelle 3.1. verdeutlicht auch die Struktur der Lösungsmenge: die Fälle pp, pm, mp und mm bilden Rechtecke und Streifen. Die Struktur entsteht durch die Positionen der Primzahlen in |H|. Die pp-Rechtecke bilden alle zusammen die Lösungsmenge. Die Träger-Darstellung beweist, dass die gezeigte Lösungsmenge vollständig ist. Die obigen drei Kongruenzen entsprechen den Kombinationen der Halbachsen $H_+$ und $H_-$: 1. Quadrant "$\equiv +2 \bmod 3\cdot2$", 3. Quadrant "$\equiv -2 \bmod 3\cdot2$", 2. und 4. Quadrant "$\equiv 0 \bmod 3\cdot2$".

<u>Die Lösungen der (2GH) zu einer bestimmten geraden Zahl</u>

Wir verwenden das Primzahl-Gitter als Koordinaten-System. Dann ist die GOLDBACH-Gleichung ‚$p+q = 2\cdot n$', allgemein $x = -y+c$, eine Geraden-Gleichung. Der geometrische Ort der Lösungen von $p+q = 2\cdot n$ für eine bestimmte gerade Zahl $2\cdot n$ ist ein Geraden-Stück; es wird gebildet von zwei Komponenten-Abschnitten der Länge m:

**Definition 3.1.** *Der Komponenten-Abschnitt $H_{...}(m)^{...}$ ist ein Anfangs-Abschnitt d.h. Abschnitt der ersten m Glieder von $H_+$ oder $H_-$. 1) Er ist ein m-tupel. 2) Suffix +: er stammt von $H_+$, Suffix -: er stammt von $H_-$. 3) Exponent +1: Reihenfolge der Glieder ‚steigt'. Exponent -1: … ‚fällt'. 4) Seine Glieder sind $h_1, h_2, h_3, \ldots, h_s, \ldots h_m$. 5) s ist Position des Gliedes innerhalb des Abschnitts bei ‚steigender' Reihenfolge der Träger. 6) Bei ‚fallender' Reihenfolge gilt für Position i und Träger s: $s + i = m + 1$.*

Wir betrachten: 1. Quadrant, $H_+ \times H_+$. Komponenten-Abschnitt $H_+(m)^{+1} = (h_1, h_2, h_3, \ldots h_m) = (7, 13, 19, \ldots, 3\cdot2\cdot m+1)$ auf der x-Achse $H_+$ und auf der y-Achse $H_+$. Punkte: $(h_1, h_m), (h_2, h_{m-1}), (h_3, h_{m-2}), \ldots, (h_m, h_1)$. Sie bilden eine Gerade, orthogonal zur Winkelhalbierenden des Quadranten, von einer Achse zur anderen verlaufend. Alle diese Punkte haben die gleiche Komponenten-Summe: $h_1 + h_m = h_2 + h_{m-1} = h_3 + h_{m-2}, \ldots h_m + h_1$. Sie stellen im Sinne von (2GH) die gerade Zahl $h_1 + h_m$ dar. $h_1 = 7 = 2\cdot3\cdot1+1$, $h_m = 2\cdot3\cdot m+1$, $h_1 + h_m = 2\cdot3\cdot(1+m) + 2 = 2\cdot3\cdot m+8$. Das bedeutet: 1) Die Zahlen $2\cdot n = 2\cdot3\cdot m+8$, $m = 1, 2, 3, \ldots$ i. e. 14, 20, 26, 32, … i. e. $(2\cdot3; 14)$ sind alle $\equiv 2 \pmod{2\cdot3}$ und sind genau die Zahlen, die oben für den 1. Quadranten erkannt waren. 2) In dieser Weise werden alle Zahlen dieses Typs dargestellt, jede genau 1 Mal - mit so vielen Lösungen auf dem Geradenstück, wie m erlaubt. GOLDBACH vermutet, dass mindestens eine von ihnen eine Lösung ist mit zwei Primzahl-Komponenten.

Wir treffen diese Verhältnisse - mutatis mutandis - auch in den anderen Quadranten an. Im 3. Quadranten befinden sich die Lösungen für die Zahlen $(2\cdot3; 10)$, im 2. und 4. Quadranten die Lösungen für $(2\cdot3; 12)$. Auf diese Weise schließt sich der geometrische Ort und bildet einen Rhombus. Genauer: Die Zahlen 10, 12, 14, 12 bilden einen ersten Rhombus, 16, 18, 20, 18 einen zweiten, u.s.w. in Schritten von $+2\cdot3$. Die Rhomben sind konzentrisch mit dem Zentrum (+1, +1), symmetrisch zu den Achsen, mit den Ecken auf den Achsen. Siehe Tabelle 3.1. Wir definieren:

Der GOLDBACH-*Lösungsrhombus* ist konzentrisch zu (+1, +1); die Seiten sind orthogonal zu den Winkelhalbierenden der Quadranten. Die Seiten werden aufgespannt von $H_+$ und $|H_-|$: 1. Rhombus-Seite.: $(H_+(m)^{+1}, H_+(m)^{-1})$, 2. Q.: $(H_-(m)^{+1}, H_+(m)^{-1})$, 3. Q.: $(H_-(m)^{+1}, H_-(m)^{-1})$, 4. Q.: $(H_+(m)^{+1}, H_-(m)^{-1})$. Die Punkte einer Seite sind symmetrische Paare und haben die gleiche Komponenten-Summe. Sie stellen eine gerade Zahl im Sinne der (2GH) dar. Die Träger laufen von 0 zur Summe der Träger. Es gibt keine andere Darstellung dieser geraden Zahl mit Primzahl-Komponenten > 3. Das Ergebnis ist damit:

**Satz 3.8.** *Der geometrische Ort der Lösungsmenge der 2-dimensionalen GOLDBACH-Vermutung (2GH) ist eine Teilmenge des linearen Systems von ‚GOLDBACH Lösungs-Rhomben':*

*1. Quadrant: $x+y = 2\cdot n \in (3\cdot2; 14) \equiv +2 \pmod{3\cdot2}$. Die (x, y) werden aufgespannt von $(H_+(m)^{+1}, H_+(m)^{-1})$, mit $2\cdot n = 3\cdot2\cdot(m-1)+14$. Sie bilden eine Rhombus-Seite. Die Komponenten $x = 3\cdot2\cdot i+1$ und $y = 3\cdot2\cdot j+1$, $i=1, 2, 3, \ldots$, $i+j = m+1$, $2\cdot n = 3\cdot2\cdot(m-1)+14$, sind symmetrisch.*

*2. Quadrant: $x+y = 2\cdot n \in (3\cdot2; 12) \equiv 0 \pmod{3\cdot2}$. Die (x, y) werden aufgespannt von $(|H_-|(m)^{+1}, H_+(H)^{-1})$, mit $2\cdot n = 12+3\cdot2\cdot(m-1)$. Sie bilden eine Rhombus-Seite. Die Komponenten $x = 3\cdot2\cdot i+1$ und $y = 3\cdot2\cdot j+1$, mit $i=1, 2, 3, \ldots$, $i+j = m+1$, $2\cdot n = 12 +3\cdot2\cdot(m-1)$, sind symmetrisch.*

*3. Quadrant: $x+y = 2\cdot n \in (3\cdot2; 10) \equiv -2 \pmod{3\cdot2}$. Die (x, y) werden aufgespannt von $(H_-(m)^{+1}, H_-(m)^{-1})$, mit $2\cdot n = 10+3\cdot2\cdot(m-1)$. Sie bilden eine Rhombus-Seite. Die Komponenten $x = 3\cdot2\cdot i+1$ und $y = 3\cdot2\cdot j+1$, mit $i=1, 2, 3, \ldots$, $i+j = m+1$, $2\cdot n = 10 +3\cdot2\cdot(m-1)$, sind symmetrisch.*

*4. Quadrant: Wie im 2. Quadranten, aber mit vertauschten Komponenten-Abschnitten.*

2. Quadrant: $E_0$ 1. Quadrant: $E_{+2}$

| | | | | | | | | | | | | | | | | | | | | | | |
|---|---|---|---|---|---|---|---|---|---|---|---|---|---|---|---|---|---|---|---|---|---|---|
| | | | | | | | | | | | ↑ | | | | | | | | | | | |
| | 120. | | 108 | 102 | . | 90 | 84 | 78 | 72 | 66 | 61 | 68 | 74 | 80 | . | 92 | 98 | 104 | . | . | 122 | |
| | | | | | | | | | | | 55 | | | | | | | | | | | |
| | | | | | | | | | | . | 49 | | | | | | | | | | | |
| | 102. | | 90 | 84 | . | 72 | 66 | 60 | 54 | 48 | 43 | 50 | 56 | 62 | . | 74 | 80 | 86 | . | . | 104 | |
| | 96 | . | 84 | 78 | . | 66 | 60 | 54 | 48 | 42 | 37 | 44 | 50 | 56 | . | 68 | 74 | 80 | . | . | 98 | |
| | 90 | . | 78 | 72 | . | 60 | 54 | 48 | 42 | 36 | 31 | 38 | 44 | 50 | . | 62 | 68 | 74 | . | . | 92 | |
| | | | | | | | | | | . | 25 | | | | | | | | | | | |
| | 78 | . | 66 | 60 | . | 48 | 42 | 36 | 30 | 24 | (19) | 26 | 32 | 38 | . | 50 | 56 | 62 | . | . | 80 | |
| | 72 | . | 60 | 54 | . | 42 | 36 | 30 | 24 | (18) | <u>13</u> | (20) | 26 | 32 | . | 44 | 50 | 56 | . | . | 74 | |
| | 66 | . | 54 | 48 | . | 36 | 30 | 24 | (18) | <u>12</u> | 7 | <u>14</u> | (20) | 25 | . | 37 | 43 | 49 | . | . | 67 | |
| ← | 59 | 55 | 47 | 41 | 35 | 29 | 23 | (17) | <u>11</u> | 5 | 1 | 7 | <u>13</u> | (19) | 25 | 31 | 37 | 43 | 49 | 55 | 61 | → |
| | 64 | . | 52 | 46 | . | 34 | 28 | 22 | (16) | <u>10</u> | 5 | <u>12</u> | (18) | 24 | . | 36 | 42 | 48 | . | . | 66 | |
| | 70 | . | 58 | 52 | . | 40 | 34 | 28 | 22 | (16) | <u>11</u> | (18) | 24 | 30 | . | 42 | 48 | 54 | . | . | 72 | |
| | 76 | . | 64 | 58 | . | 46 | 40 | 34 | 28 | 22 | (17) | 24 | 30 | 36 | . | 48 | 54 | 60 | . | . | 78 | |
| | 82 | . | 70 | 64 | . | 52 | 46 | 40 | 34 | 28 | 23 | 30 | 36 | 42 | . | 54 | 60 | 66 | . | . | 84 | |
| | 88 | . | 76 | 70 | . | 58 | 52 | 46 | 40 | 34 | 29 | 36 | 42 | 48 | . | 60 | 66 | 72 | . | . | 90 | |
| | | | | | | | | | | . | 35 | | | | | | | | | | | |
| | 100 | . | 88 | 82 | | 70 | 64 | 58 | 52 | 46 | 41 | 48 | 54 | 60 | . | 72 | 78 | 84 | . | . | 102 | |
| | 106 | . | 94 | 88 | . | 76 | 70 | 64 | 58 | 52 | 47 | 54 | 60 | 66 | . | 78 | 84 | 90 | . | . | 108 | |
| | 112 | . | 100 | 94 | . | 82 | 76 | 70 | 64 | 58 | 53 | 60 | 66 | 72 | . | 84 | 90 | 96 | . | . | 114 | |
| | 118 | . | 106 | 100 | . | 88 | 82 | 76 | 70 | 64 | 59 | 66 | 72 | 78 | . | 90 | 96 | 102 | . | . | 120 | |
| | | | | | | | | | | | ↓ | | | | | | | | | | | |

3. Quadrant: $E_{-2}$ 4. Quadrant: $E_0$

(2GH): Lösungsmenge und Struktur der Lösungsmenge der (2GH) in Original-Darstellung
zwei Beispiele markiert: <u>10</u>, <u>12</u>, <u>14</u>, <u>12</u> und (16), (18), (20), (18)
Tabelle 3.1.

Die Seiten eines GOLDBACH Lösungsrhombus verlaufen ‚diagonal'. Deshalb sind sie in der Lage, im Primzahl-Gitter über die Streifen zu springen, die von den Vielfachen gebildet werden. Dieses Phänomen zeigt, dass es in der GOLDBACH-Diskussion vorteilhaft ist, die Positionen der Primzahlen und Vielfachen zu bewahren. Ohne Bewahrung der Positionen würde die Struktur nicht zum Vorschein kommen. - Außerdem liefert dieses Phänomen den Grund, weshalb bei wachsendem m die Anzahl Primzahlen und die Anzahl Vielfache nicht streng monoton wachsen, sondern das nur mit gewissen Schwankungen tun.

## § 4. Geometrischer Ort der Lösungsmenge der Primzahlzwillinge. Verallgemeinerung.

Primzahlzwillinge sind Primzahlen p und q, deren Differenz q - p = 2 ist. Wenn (q, p) ein Primzahlzwilling ist und q < p, so ist q Glied von $H_-$ und p Glied von $H_+$. Wir nennen das Problem der Primzahlzwillinge (PT). Bei Verwendung von H und Primzahlgitter ist die definierende Gleichung nicht p - q = 2, sondern p + q = 2, mit negativem q. - Versteht man p + q = 2 als Geraden-Gleichung, dann verläuft die Gerade $p = (-1)\cdot q + 2$ mit der Steigung tg $\alpha = -1$ durch den Punkt (+1, +1). Jeder Punkt dieser Geraden hat die Koordinaten-Summe x + y = 2. Die Gerade ist die Winkelhalbierende des 2. und 4. Quadranten. Die Gerade enthält alle Zwillinge, insbesondere die Primzahlzwillinge. Sie enthält alle Arten von Kombinationen pp, pm, mp, mm. (Siehe Tab. 4.1.) Also gilt:

**Satz 4.1.** *Im Primzahlgitter $H \times H$ ist der geometrische Ort der Lösungsmenge der Primzahlzwillinge (PT) die Winkelhalbierende des 2. (bzw. 4.) Quadranten: die Gerade $p = (-1)\cdot q + 2$.*

Verallgemeinerung:

Man kann den Begriff 'Primzahlzwilling' verallgemeinern. Hierzu habe die Differenz p - q andere Werte als 2. Das Primzahl-Gitter zeigt, dass diese Werte $p + q = 3\cdot 2\cdot i + 2$, $i \in N$, sein müssen. Dann sind die geometrischen Örter der Lösungsmengen Geraden mit den Eigenschaften: 1) sie sind Parallelen zur Winkelhalbierenden des 2. Quadranten, 2) sie schneiden diejenigen Achsenpunkte auf der y-Achse, die Glieder von H sind.

Man kann die Verallgemeinerung der Zwillinge auch in den 1. Quadranten verlegen. Dort sind p und q positiv. Man setzt an: p - q = 0 bzw. allgemein $p - q = 3\cdot 2\cdot i$, $i \in N$. Dann erhält man zunächst die Winkelhalbierende des 1. Quadranten. In der Verallgemeinerung ergibt sich die Parallelenschar mit den Eigenschaften: 1) sie sind Parallelen zur Winkelhalbierenden des 1. Quadranten, 2) sie schneiden diejenigen Achsenpunkte auf der y-Achse, die Glieder von H sind, 3) sie bedecken schließlich alle vier Quadranten. - Wir nennen dieses Problem (gPT), generalisierte Primzahl-Zwillinge. Siehe Tabelle 4.1. Es gilt:

**Satz 4.2.** *Der geometrische Ort der Lösungen von (gPT) im Primzahl-Gitter: Die Parallelenscharen $p = (-1)\cdot q + 3\cdot 2\cdot i + 2$, i fest $\in N$ bzw. $p = q + 3\cdot 2\cdot i$, i fest $\in N$, mit: 1) verlaufen parallel zur Winkelhalbierenden des 2. Quadranten, 2) bedecken schließlich alle vier Quadranten, 3) schneiden (auf der y-Achse) die Glieder von H.*

| | | | | | | | | | | | | | | | | | | | | | | | |
|---|---|---|---|---|---|---|---|---|---|---|---|---|---|---|---|---|---|---|---|---|---|---|---|
| | - | - | - | - | - | - | - | - | - | - | - | - | - | - | - | - | - | - | - | - | - | 115 | - |
| | - | -4 | <u>2</u> | 8 | - | 20 | 26 | - | 38 | - | 50 | 56 | 62 | 68 | - | 80 | 86 | 92 | 98 | 104 | - | 109 | p |
| | - | -10 | -4 | <u>2</u> | - | 14 | 20 | - | 32 | - | 44 | 50 | 56 | 62 | - | 74 | 80 | 86 | 92 | 98 | - | 103 | p |
| | - | -16 | -10 | -4 | - | 8 | 14 | - | 26 | - | 38 | 44 | 50 | 56 | - | 68 | 74 | 80 | 86 | 92 | - | 97 | p |
| | - | - | - | - | - | - | - | - | - | - | - | - | - | - | - | - | - | - | - | - | - | 91 | - |
| | - | - | - | - | - | - | - | - | - | - | - | - | - | - | - | - | - | - | - | - | - | 85 | - |
| | - | -34 | -28 | -22 | - | -10 | -4 | - | 8 | - | 20 | 26 | 32 | 38 | - | 50 | 56 | 62 | 68 | 74 | - | 79 | p |
| | - | -40 | -34 | -28 | - | -16 | -10 | - | <u>2</u> | - | 14 | 20 | 26 | 32 | - | 44 | 50 | 56 | 62 | 68 | - | 73 | p |
| | - | -46 | -40 | -34 | - | -22 | -16 | - | -4 | - | 8 | 14 | 20 | 26 | - | 38 | 44 | 50 | 56 | 62 | - | 67 | p |
| | - | -52 | -46 | -40 | - | -28 | -22 | - | -10 | - | <u>2</u> | 8 | 14 | 20 | - | 32 | 38 | 44 | 50 | 56 | - | 61 | p |
| | - | - | - | - | - | - | - | - | - | - | - | - | - | - | - | - | - | - | - | - | - | 55 | - |
| | - | - | - | - | - | - | - | - | - | - | - | - | - | - | - | - | - | - | - | - | - | 49 | - |
| | - | -70 | -64 | -58 | - | -46 | -40 | - | -28 | - | -16 | -10 | -4 | <u>2</u> | - | 14 | 20 | 26 | 32 | 38 | - | 43 | p |
| | - | -76 | -70 | -64 | - | -52 | -46 | - | -34 | - | -22 | -16 | -10 | -4 | - | 8 | 14 | 20 | 26 | 32 | - | 37 | p |
| | - | -82 | -76 | -70 | - | -58 | -52 | - | -40 | - | -28 | -22 | -16 | -20 | - | <u>2</u> | 8 | 1 4 | 20 | 26 | - | 31 | p |
| | - | - | - | - | - | - | - | - | - | - | - | - | - | - | - | - | - | - | - | - | - | 25 | - |
| | - | -94 | -88 | -82 | - | -70 | -64 | - | -52 | - | -40 | -34 | -28 | -22 | - | -10 | -4 | <u>2</u> | 8 | 14 | - | 19 | p |
| | - | -100 | -94 | -88 | - | -76 | -70 | - | -58 | - | -46 | -40 | -34 | -28 | - | -16 | -10 | -4 | <u>2</u> | 8 | - | 13 | p |
| | - | -106 | -100 | -94 | - | -82 | -76 | - | -64 | - | -52 | -46 | -40 | -34 | - | -22 | -16 | -10 | -4 | <u>2</u> | - | 7 | p |
| | - | - | - | - | - | - | - | - | - | - | - | - | - | - | - | - | - | - | - | - | - | 1 | - |
| H- ← | -119 | -113 | -107 | -101 | -95 | -89 | -83 | -77 | -71 | -65 | -59 | -53 | -47 | -41 | -35 | -29 | -23 | -17 | -11 | -5 | 1 | | |
| | - | p | p | p | - | p | p | - | p | - | p | p | p | p | - | p | p | p | p | p | | | |

2. Quadrant
(PT): Primzahl-Gitter H×H, Punkte x + y, nur pp-Punkte sind dargestellt.
Für Primzahlzwillinge gilt p + q = 2 (Unterstreichung); sie bilden eine Teilmenge der Winkelhalbierenden
Tabelle 4.1.

**§ 5. Geometrischer Ort der Lösungsmenge der PRACHAR Primzahl-zwillinge. Verallgemeinerung.**

Im Umgang mit den PRACHAR-Primzahlzwillingen, (PPT) genannt, verwenden wir H zunächst in der Form H = (-1)·(± 3·2; -1). Dabei ist das Primzahl-Gitter (-1)·H × (-1)·H = ((-1)·(± 3·2; -1)) × ((-1) · (± 3·2; -1)). Das ist das bereits erwähnte Primzahl-Gitter mit dem Ursprung (-1, -1). Wir suchen den geometrischen Ort der Lösungsmenge von (PPT) im Primzahl-Gitter (-1)·H × (-1)·H. Die von PRACHAR gegebene Definition von (PPT) lautet p = ½·(q - 1). Beispiele: ...(+23, +47), (+17, +35), (+11, +23), (+5, +11), (-1, -1), (-7, -13), (-13, -25), (-19, -37), (-25, -49), (-31, -61), (-37, -73), (-43, -85), ... Wir verstehen diese Kombinationen von Komponenten als Punkte des Primzahl-Gitters. D. h.: der geometrische Ort der Lösungsmenge von (PPT) ist die Gerade p = ½·(q - 1) im Primzahl-Gitter. Sie enthält die Punkte (-1, -1) und (5, 11) und hat die Steigung tg $\alpha$ = ½. Sie enthält Punkte mit zwei Primzahl-Komponenten; siehe Tabelle 5.1.

Harmonisierung: Soll (PPT) im Primzahl-Gitter H × H erscheinen, ändert sich die originäre Definitionsgleichung von p = ½·(q - 1) zu p = ½·(q + 1). Abgesehen von Vorzeichen ändert sich die Lösungsmenge nicht.

| | | | | | | | | | | | | | | | | | | | |
|---|---|---|---|---|---|---|---|---|---|---|---|---|---|---|---|---|---|---|---|
| -55 | . | . | . | mm | . | . | . | . | . | mm | . | . | . | mm | . | . | . | mm | mm |
| -49 | . | . | . | mm | . | . | . | . | . | mm | . | . | . | mm | . | . | . | mm | mm |
| -43 | pp | pp | pp | . | pp | pp | pp | pp | pp | . | pp | pp | pp | . | pp | pp | pp | | |
| -37 | pp | pp | pp | . | pp | pp | pp | pp | pp | . | pp | pp | • | . | pp | pp | pp | | |
| -31 | pp | pp | pp | . | pp | pp | pp | pp | pp | . | pp | pp | pp | . | pp | pp | pp | | |
| -25 | . | . | · | mm | . | . | . | . | . | mm | . | . | . | mm | . | . | . | .mm | mm |
| -19 | pp | pp | pp | . | pp | pp | pp | pp | pp | . | pp | pp | pp | . | pp | pp | pp | | |
| -13 | pp | pp | pp | . | pp | pp | pp | pp | pp | . | • | pp | pp | . | pp | pp | pp | | |
| - 7 | pp | pp | pp | . | pp | pp | pp | pp | pp | . | pp | pp | pp | . | pp | pp | pp | | |
| -1 | . | . | . | mm | . | . | . | . | . | . | . | . | . | mm | . | . | . | mm | mm |
| 5 | pp | pp | pp | . | pp | pp | pp | pp | pp | . | pp | pp | pp | . | pp | pp | pp | | |
| 11 | pp | pp | pp | . | pp | pp | pp | pp | • | . | pp | pp | pp | . | pp | pp | pp | | |
| 17 | pp | pp | pp | . | pp | pp | pp | pp | pp | . | pp | pp | pp | . | pp | pp | pp | | |
| 23 | pp | pp | pp | . | pp | pp | pp | • | pp | . | pp | pp | pp | . | pp | pp | pp | | |
| 29 | pp | pp | pp | . | pp | pp | pp | pp | pp | . | pp | pp | pp | . | pp | pp | pp | | |
| 35 | . | . | . | mm | . | . | . | . | . | mm | . | . | . | mm | . | . | . | mm | mm |
| 41 | pp | pp | pp | . | pp | pp | pp | pp | pp | . | pp | pp | pp | . | pp | pp | pp | | |
| 47 | pp | pp | pp | . | pp | • | pp | pp | pp | . | pp | pp | pp | . | pp | pp | pp | | |
| 53 | pp | pp | pp | . | pp | pp | pp | pp | pp | . | pp | pp | pp | . | pp | pp | pp | | |
| | | | | | | | | | | ↑ | | | | | | | | | |
| (-1)·H → | 53 | 47 | 41 | 35 | ·29 | 23 | 17 | 11 | 5 | -1 | -7 | -13 | -19 | -25 | -31 | -37 | -43 | -49 | -55 |

(PPT): Primzahl-Gitter (-1) ·H × (-1) ·H, pp-Punkte,
Lösungen von (PPT) bei •
Tabelle 5.1.

Verallgemeinerung: Der (hier eingeführte) Begriff der 'PRACHAR-Primzahlzwillinge' kann verallgemeinert werden. Dazu arbeiten wir in H × H. Die Gitter-Punkte sind dann mit Paaren (a, b) besetzt, wobei a und b Glieder von H sind. a nennt die Glieder der x-Achse, b die der y-Achse. Sind a und b Primzahlen, so handelt es sich um eine Lösung. - Die von PRACHAR gegebene Gleichung lautet $y = ½·(x + 1)$. Sie enthält im Primzahl-Gitter die Punkte (+1, +1) und (+13, +7) und ist durch sie festgelegt. Eine Liste von Beispielspunkten lässt sich beliebig erweitern: ... (-35, -17), (-23, -11), (-11, -5), (+1, +1), (+13, +7), (+25, +13), (+37, +19), ... - Zur Verallgemeinerung drehen wir die Gerade um den Punkt (+1, +1) so, dass sie wieder durch einen zweiten Gitterpunkt läuft. Eine Gerade durch zwei gegebene Punkte $(x_1, y_1)$ und $(x_2, y_2)$ folgt der *2-Punkt-Gleichung*: $(y - y_1)/(x - x_1) = (y_2 - y_1)/(x_2 - x_1)$. Bei der Drehung halten wir einen der beiden Punkte fest, und zwar $x_1 = +1$, $y_1 = +1$. Dann lautet die Geraden-Gleichung:
$(y - 1)/(x - 1) = (y_2 - 1)/(x_2 - 1)$.
Wählt man $x_2 = +13$, $y_2 = +7$, so erhält man die von PRACHAR vorgegebene Gerade: $(y - 1)/(x - 1) = (7 - 1)/(13 - 1)$, d.h. $(y - 1) = ½·(x - 1)$, d.h. $y = ½·(x + 1)$. Sie enthält z. B. die o. a. Punkte. Sie verläuft im 1. und 3. Quadranten. Offenbar enthält sie Lösungen, z. B.: (-23, -11) und (-11, -5) und (+13, +7) und (+37, +19) und ...
Wählt man $x_2 = +7$, $y_2 = +7$, so erhält man die Winkelhalbierende x = y. Sie enthält trivialerweise ∞ viele Lösungen, nämlich für jede Primzahl genau eine. - Wählt man als Erst-Punkt (+1, +1) und als Zweit-Punkte Punkte mit $x_2 > y_2$, so erhält man alle Geraden, die im 1. und 3. Quadranten verlaufen, (+1, +1) enthalten, den Zweit-Punkt $(x_2, y_2)$ enthalten und flacher verlaufen als die Winkelhalbierende. Wählt man Zweit-Punkte mit $x_2 < y_2$, so erhält man die Geraden, die steiler verlaufen. - Die gleiche Prozedur kann man auch im 2. und somit im 4. Quadranten ansetzen. Z.B. sei etwa $(x_2, y_2) = (-11, +7)$. Dann lautet die Geraden-Gleichung $(y - 1)/(x - 1) = (+7 - 1)/(-11 - 1) = -½$, d.h.
$y = ½·(-x +3)$. Beispielspunkte: ..., (-35, +19), (-23, +13), (-11, +7), (+1, +1), (+13, -5), (+25, -11), (+37, -17), ... Man erhält schließlich ein Geraden-Bündel mit Schnittpunkt (+1, +1) und den Gitter-Punkten als Zweit-Punkten. Ausgenommen sind die Achsen, da 1 nicht Primzahl ist - oder: 1 spielt die Rolle einer Primzahl. (Hier zeigt sich erneut: die Szenerie der Primzahlen ist besser abgerundet, wenn auch 1 eine Primzahl ist.) - Die gleiche Prozedur kann man an j e d e m Gitter-Punkt veranstalten.

**Satz 5.1.** *Der geometrische Ort der Lösungsmenge der PRACHAR-Primzahlzwillinge (PPT) ist die Gerade y = ½·(x + 1), d. h. durch den Punkt (+1, +1) mit der Steigung ½. - Der geometrische Ort der Lösungsmenge der (gPPT) ist das Geraden-Bündel mit Schnitt-Punkt (+1, +1) und den Gitter-Punkten als Zweit-Punkten. - Der geometrische Ort der Lösungsmenge der total verallgemeinerten (gPPT) sind die vergleichbaren Geraden-Bündel in jedem Gitter-Punkt.*

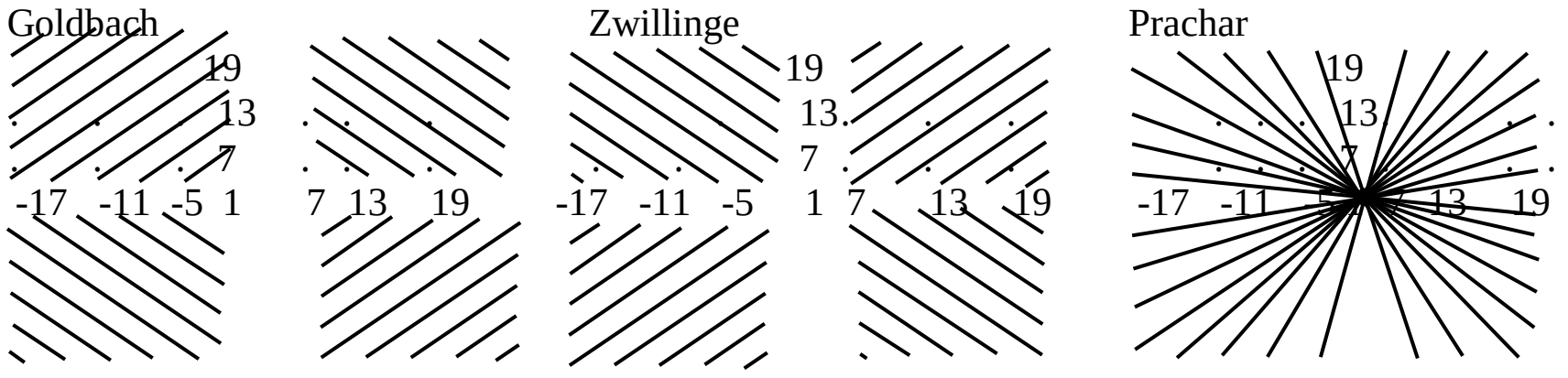

Skizze: die geometrischen Örter der Lösungsmengen im Primzahl-Gitter (verallgemeinert)
Tabelle 5.2.

**§ 6. Geometrischer Ort der Lösungsmenge der 3-gliedrigen GOLDBACH-Vermutung**.
<u>Szenerie</u>: Wir nennen das Problem (3GH). Wir schließen wieder 2 und 3 von den Lösungskomponenten aus. Wir verstehen das Problem als ein 3-dimensionales. Wir bilden das 3-dimensionale Primzahl-Gitter H×H×H. Da (3GH) alle Primzahlen als positiv voraus setzt, sprechen wir nicht von H sondern von |H|. Um die Betragstriche zu vermeiden, sei vereinbart, dass in diesem Kapitel immer |H| gemeint ist, wenn H geschrieben wird. - Die Achsen seien zueinander orthogonal mit äquidistanter Einteilung. Wir nennen sie x-, y- und z-Achse. Ursprung ist der Punkt (+1, +1, +1). - Alle Komponenten haben als Glieder von H die Form (3·2·s ± 1). Suffixe am Träger geben die Herkunftsachse an. Ein Punkt hat dann die Koordinaten $((3·2·s_x ± 1), (3·2·s_y ± 1), (3·2·s_z ± 1))$. Wir besetzen die Gitter-Punkte mit der Koordinaten-Summe.
<u>Lösungs-Kubus</u>: H ist die Vereinigung von $H_-$ und $H_+$. Also ist der Koordinaten-Kubus in 8 Teil-Kuben eingeteilt, aufgespannt von $H_-$ und $H_+$: 'Oberer Halb-Kubus', beginnend in $H_+ × H_+ × H_+$, Zählung der Teil-Kuben entgegen dem Uhrzeigersinn: $H_+ × H_+ × H_+$ ist Teil-Kubus 1. $H_- × H_+ × H_+$, ist Teilkubus 2. $H_- × H_- × H$ ,ist Teilkubus 3, $H_+ × H_- × H_+$ ist Teilkubus 4.'Unterer Halb-Kubus', beginnend in Teil-Kubus 5, diametral zu Teil-Kubus 1, Zählung entgegen dem Uhrzeigersinn: $H_- \times H_- \times H_-$ ist Teil-Kubus 5 diametral zu 1, $H_+ \times H_- \times H_-$ ist Teilkubus 6, $H_+ \times H_+ \times H_-$ ist Teilkubus 7, $H_- \times H_+ \times H_-$ ist Teilkubus 8.

<u>Linearität</u>: Wir diskutieren im 1. Teil-Kubus die Gleichung $x + y + z = C_1$. Die Gitter-Punkte, die diese Gleichung erfüllen, liegen auf einer Ebene. Dabei geht es um denjenigen Teil der Ebene, der im 1.Teil-Kubus liegt. Er wird aufgespannt von den Eck-Punkten $(C_1, 0, 0)$, $(0, C_1, 0)$, $(0, 0, C_1)$: er ist ein Dreieck, dessen Seiten die Eckpunkte verbinden. Mit $x = (3\cdot2\cdot s_1+1)$, $y = (3\cdot2\cdot s_2+1)$, $z = (3\cdot2\cdot s_3+1)$ lautet die Ebenen-Gleichung $(3\cdot2\cdot s_1+1) + (3\cdot2\cdot s_2+1) + (3\cdot2\cdot s_3+1) = C_1$.

2. Teil-Kubus
$u = 25 = 3\cdot2\cdot4 + 1$
$\Sigma s_i = s_x+s_y+s_z$ = constans = 4

$H_+$

1. Teil-Kubus
$u = 27 = 3\cdot2\cdot4 + 3$
$\Sigma s_i=s_x+s_y+s_z$=constans=4

(1,25,1)
*0 4 0*

(5,9,1) ¦ (1,19,7) (7,19,1)
*1, 3, 0* ¦ *0 3 1* *1 3 0*

(11,13,1) (5,13,7) ¦ (1,13,13) (7,13,7) (13,13,1)
*2 2 0* *1 2 1* ¦ *0 2 2* *1 2 1* *2 2 0*

(17,7,1) (11,7,7) (5,7,13) ¦ (1,7,19) (7,7,13) (13,7,7) (19,7,1)
*3 1 0* *2 1 1* *1 1 2* ¦ *0 1 3* *1 1 2* *2 1 1* *3 1 0*

(23,1,1) (17,1,7) (11,1,13) (5,1,19) ¦ (1,1,25) (7,1,19) (13,1,13) (19,1,7) *0 0*
*4 0 0* *3 0 1* *2 0 2* *1 0 3* ¦ *0 0 4* *1 0 3* *2 0 2* *3 0 1*

$H_-$ - - - - - - - - - - - - - - - - - - - - - - - - - - - - - - - - - - - - - - - - - - - - - - - - - - - - $H_+$

(17,5,1) (11,5,7) (5,5,13) ¦ (1,5,19) (7,5,13) (13,5,7) ( 19,5,1)
*3 1 0* *2 1 1* *1 1 2* ¦ *0 1 3* *1 1 2* *2 1 1* *3 1 0*

(11,11,1) (5,11,7) ¦(1,11,13) (7,11,7) (13,11,1)
*2 2 0* *1 2 1* ¦ *0 2 2* *1 2 1* *2 2 0*

(5,17,1) ¦ (1,17,7) (7,17,1)
2 2 0 ¦ 0 2 2 2 2 0

¦(1,23,1)
¦*0 4 0*

3. Teil-Kubus
$u = 23 = 3\cdot2\cdot4 - 1$
$\Sigma s_i = s_x+s_y+s_z$=constans=4

$H_-$

4. Teil-Kubus
$u = 25 = 3\cdot2\cdot4 + 1$
$\Sigma s_i = s_x+s_y+s_z$=constans=4

(3GH): Lösungen der 3-gliedrigen GOLDBACHschen Vermutung.
Oktaeder, 'obere Hälfte', Sicht von 'oben', 'rechts' abgeschnitten.
Beispiel: Träger-Summe $\Sigma s_i = s_x + s_y + s_z$ = constans = 4.
Tabelle 6.1.

Bei dem geschilderten Ebenen-Stück ist auch die Summe der Träger konstant: $s_1 + s_2 + s_3 = s$ ist konstant. Für die Komponenten-Summe gilt im 1. Teil-Kubus $(3\cdot2\cdot s_1+1) + (3\cdot2\cdot s_2+1) + (3\cdot2\cdot s_3+1) = (3\cdot2\cdot s_1 + (3\cdot2\cdot s_2 + 3\cdot2\cdot s_3) +1+1+1 = 3\cdot2\cdot(s_1 + s_2 + s_3) +1+1+1$. Man kann also den Summanden $s_1+s_2+s_3$ festhalten und den Summanden 1+1+1 variieren.

<u>Oktaeder</u>: Wir halten die Träger-Summe $s_1+s_2+s_3$ fest und variieren die Reste-Summe im Sinne von $\pm1$, $\pm1$, $\pm1$, wie es die Teil-Kuben vorschreiben. D.h.: wir erhalten bei fester Träger-Summe nur mit Variation der Reste 8 Ebenen-Stücke, in jedem Teil-Kubus eines. Sie bilden zusammen ein Oktaeder. Fälle: 1. Teil-Kubus: +1+1+1 = +3, 2. Teil-Kubus: -1+1+1 = +1, 3. Teil-Kubus: -1-1+1 = -1, 4. Teil-Kubus: +1-1+1 = +1, 5. Teil-Kubus: -1-1-1 = -3, 6. Teil-Kubus: +1-1-1= -1, 7. Teil-Kubus: +1+1-1= +1, 8. Teil-Kubus: -1+1-1 = -1.

Zusammenhang zwischen den ungeraden Zahlen und den Oktaedern: n gegeben, n sei die Summe der Träger der Punkte in $H\times H\times H$. Das führt zu 4 ungeraden Zahlen: $3\cdot2\cdot n + 3$, $3\cdot2\cdot n + 1$, $3\cdot2\cdot n - 1$ und $3\cdot2\cdot n - 3$ mit $n = s_1 + s_2 + s_3$; sie bilden ein Oktaeder.

Vollständigkeit: Das System der Oktaeder bedient vollständig alle Fälle, wie von (3GH) verlangt. Die 4-Tupel für n überlappen. Das sieht man, wenn man die Träger-Summe um 1 erhöht. Um die ungeraden Zahlen ganz abzudecken, ist die Träger-Summe schrittweise um 1 zu erhöhen. Dabei entsteht ein lineares System von Oktaedern der geschilderten Art. Dieses System ist der geometrische Ort der Lösungsmenge der (3GH).

<u>Lösungen</u>: Punkte, die als Komponente ein Vielfaches besitzen, sind nicht Lösung. Punkte mit Komponente 1 (gleichbedeutend mit Träger 0) sind Rand-Punkte und sind nicht Lösung. Punkte mit drei Primzahl-Komponenten sind Lösung.

Zuordnung der Ränder: Die Halb-Achse $H_-$ enthält nicht das Glied -1. Hieraus folgt die Klassifikation der Ecken der Kuben: Die Kuben sind zueinander fremd. Kuben, die den Punkt (1, 1, 1) enthalten, sind 'abgeschlossen', die anderen sind 'offen'. In diesem Sinne sind die Kuben nicht gleich. Siehe Tabelle 6.1. Beispiel: Es sei n = 4. Die darzustellenden ungeraden Zahlen sind 27 (1. Kubus), 25 (2., 4. und 7. Kubus), 23 (3., 6. und 8. Kubus) und 21 (5. Kubus). Diese 4 ungeraden Zahlen besetzen genau ein Oktaeder. Lösungen sind: Punkt (7, 13, 7) mit Träger (1, 2, 1), Punkt (7, 7, 13) mit Träger (1, 1, 2) und Punkt (13, 7, 7) mit Träger (2, 1, 1). Die Parameter der ersten Oktaeder sind: <u>$(3\cdot2\cdot n + 3)$, $(3\cdot2\cdot n + 1)$, $(3\cdot2\cdot n -1)$, $(3\cdot2\cdot n -3)$:</u> n=1: kein Okt. n=2: <u>(9), (11), (13), (15).</u> n=3: <u>(15), (17), (19),</u>

(21). n=4: (21), (23), (25), (27). n=5: (27), (29), (31), (33). n=6: (33), (35), (37), (39). ... Die 'untere' Hälfte ergibt ein analoges Bild. Es gilt:

**Satz 6.1.** *Der geometrische Ort der Lösungsmenge der (3GH) ist ein lineares System von Oktaedern. Die Oktaeder sind parallel, äquidistant mit Abstand 3·2 auf den Achsen, konzentrisch zu (1, 1, 1), die Ecken auf den Achsen. Parameter des Systems ist die Summe der Träger s; s ist fest für den einzelnen Oktaeder. s generiert 4 ungerade Zahlen 3·2·s-3, 3·2·s -1, 3·2·s +1 und 3·2·s +3, 2 < s $\in$ N, die vom s-ten Oktaeder dargestellt werden. Die Komponenten der Punkte sind Glieder von |H|. Die Bestimmung der Komponenten aus |H| geschieht durch den Oktaeder* $(3\cdot 2\cdot s_x \pm 1) + (3\cdot 2\cdot s_y \pm 1) + (3\cdot 2\cdot s_z \pm 1)$ *mit* $s = s_x + s_y + s_z$ *und* $0 \leq s_x, s_y, s_z \leq s$.

## § 7. Naheliegende Vermutungen

1. Verallgemeinerung der Primzahlzwillinge: *Im Primzahl-Gitter H×H ist der geometrische Ort der verallgemeinerten Primzahlzwillinge (PT) eine Parallelen-Schar parallel zur Winkelhalbierenden des 2. (bzw. 1.) und 4. (bzw. 3.) Quadranten. Die Parallelen verlaufen im Abstand der Gitter-Punkte auf den Achsen. Vermutung: Jede dieser Parallelen enthält unendlich viele verallgemeinerte Primzahlzwillinge.*

2. Verallgemeinerte PRACHAR Primzahlzwillinge: *Im Primzahl-Gitter H×H ist der geometrische Ort der verallgemeinerten PRACHAR Primzahlzwillinge (PPT) ein Geraden-Bündel durch den Punkt (1, 1); die Geraden haben rationale Steigungen. Vermutung: Jede dieser Geraden enthält unendlich viele verallgemeinerte PRACHAR Primzahlzwillinge. Vermutung: Das gilt in jedem Gitter-Punkt.*

3. Verallgemeinerung der Primzahlzwillinge in den Raum, *Feststellung: Man kann die in Punkt 2 und Punkt 3 angesprochenen H×H-Systeme von Geraden der geometrischen Örter der Primzahlzwillinge und PRACHAR Primzahlzwillinge in den Raum H×H×H verallgemeinern.*

4. Totale Verallgemeinerung, *Vermutung: Jede nicht achsen-parallele Gerade im Primzahl-Gitter H×H, die zwei Gitter-Punkte enthält, enthält unendlich viele Punkte mit Primzahl-Komponenten. Jede nicht-achsenparallele Ebene im Primzahl-Gitter H×H×H, die drei Gitter-Punkte enthält, enthält unendlich viele Punkte mit Primzahl-Komponenten.*

5. Definition von Grundbegriffen der Geometrie allein durch Primzahlen: *Vermutung: Man kann Grundbegriffe der Geometrie wie Gerade, Parallele usw. allein mit Primzahlen definieren. Die einzigen Voraussetzungen sind die Fragestellung von (2GH) und die Beschränkung auf homöomorphe Gegenstände.*

## § 8. Relevanz-Intervalle

In diesem § ist p > 3.

Zugehörigkeit der Primzahl-Quadrate zu $H_+$:

Nach Satz 1.4. hat p die Gestalt 3·2·k +1 oder 3·2·k - 1. Also hat $p^2$ in beiden Fällen die Gestalt (…) + 1. D.h. $p^2$ ist Glied von $H_+$:

**Satz 8.1.** *Die Primzahl-Quadrate liegen in* $H_+$.

**Satz 8.2.** $p^{2n} \in H_+$, $p^{2n+1} \in H_-$.

Intervall-Einteilung mittels der Primzahl-Quadrate:

Sei p, q $\in$ P und p < q und p, q aufeinander folgend in P. Wir definieren und betrachten in P das Intervall zwischen $p^2$ u. $q^2$, $p^2$ eingeschlossen, $q^2$ ausgeschlossen, und nennen es ‚Relevanz-Intervall' oder ‚$p^2$-Intervall':

**Def. 8.1.** *Relevanz-Intervall, syn.* $p^2$*-Intervall: p, q* $\in$ *P und p < q und p, q aufeinander folgend in P:*

*Relevanz-Intervall* $(H_+; p^2)$ *in* $H_+ := \{ x \mid (x \in H_+) \wedge (p^2 \leq x < q^2)\}$.

*Relevanz-Intervall* $(H_-; p^2)$ *in* $H_- := \{ y \mid (y \in H_-) \wedge\ y = x - 2 \wedge\ x \in (H_+; p^2)\}$.

Definition der im $p^2$-Intervall 'relevanten Primzahlen', das sind die $p_i$ zwischen 3 und p:

**Def. 8.2.** *Die relevanten Primzahlen* $p_i$ *in* $(H_+; p^2)$ *u.* $(H_-; p^2)$ *sind die Primzahlen* $p_i$ *mit* $3 < p_i \leq p$.

Die relev. $p_i$ heißen 'relev.', weil ein Vielf. immer Vielf. von einer relev. $p_i$ ist. Die anderen sind sozusagen irrelevant:

**Satz 8.3.** *Es ist hinreichend, in* $(H_+; p^2)$ *u.* $(H_-; p^2)$ *Vielfache als Vielfache der dortigen relevanten* $p_i$ *zu behandeln.*

Das bedeutet: in den Vielfachen von $(H_+; p^2)$ und $(H_-; p^2)$ sind nur die $p_i$ relevant:

Differenz q – p und Länge L des Intervalls:

**Def. 8.3.** *Abstand f, Länge L: Abstand f := q - p. f ist gerade. L bezeichnet die Länge (Anzahl Glieder) des Intervalls ($H_+$; $p^2$). L := ($q^2$ - $p^2$)/3·2. - Seien $p_1 < p_2$ Primzahlzwillinge. Dann gilt: 1.) $p_1 \in H_-$ u. $p_2 \in H_+$, 2.) Das Glied zwischen $p_1$ u. $p_2$ wird von 3 geteilt, 3.) Nur bei Primzahlzwillingen ist L < p; sonst ist p < L.*

Rest-Systeme in Z, N, H, $H_+$ u. $H_-$:
**Satz 8.4.** *p aufeinander folgende Glieder von Z, N, H, $H_+$, $H_-$ lassen ein vollständiges Restsystem modulo p.*
**Satz 8.5.** *In $H_-$, $H_+$ u. H wachsen die Reste (mod p) um 3·2 pro Glied. Die Summe zweier Reste, die symmetrisch zum 0-Rest liegen, ist $\equiv 0$ (mod p). $p < q \Rightarrow q$ lässt mindestens zwei Reste mehr als p.*

Die Relevanz-Intervalle prägen die Struktur von P. Rel.-Int- und Fundamentalsatz:
**Satz 8.6.** *1.) Die Vereinigung der Relevanz-Intervalle ist IHI - abgesehen von den Gliedern -1, 2 und 3. 2.) ($H_+$; $p^2$) u. ($H_-$; $p^2$) besitzen die gleiche Länge L. 3.) ($H_+$; $p^2$) enthält mindestens ein Vielfaches, u. zwar $p^2$. 4.) ($H_+$; $p^2$) u. ($H_-$; $p^2$) haben dieselben relevanten Primfaktoren. 5.) Die Primfaktor-Zerlegung eines Vielfachen im Rel.-Int. enthält mindestens einen relevanten Primfaktor. 6.) Die Vielfachen im Rel.-Int. sind Vielfache von welchen der relevanten Primfaktoren. 7.) In jedem weiteren Rel.-Int. tritt genau ein relevanter Primfaktor hinzu, u. zwar p.*

Jede p generiert eine arithmetische Folgen, aZ ($\pm$3·2·p; $p^2$):
Die in H enthaltene Menge der Vielfachen ist wie ein ‚Zopf' aus arithmetischen Folgen geflochten. Jedes p generiert die arithmetische Folge der Vielfachen von p. H selbst ist eine arithmetische 2-Ast-Folge: in der Gestalt H = ($\pm$3·2; +1) beginnt sie mit dem Glied +1. Anfangsglied darf jedes ihrer Glieder sein. Also auch p, also auch $p^2$: H = ($\pm$3·2; $p^2$). Dann ist ($\pm$3·2·p; $p^2$) die in ($H_+$; $p^2$) mit $p^2$ beginnende arithmetische Folge der Vielfachen von p, sofern sie in H liegen. Per Definition enthält ($H_+$; $p^2$) genau ein Primzahl-Quadrat: das Anfangsglied $p^2$. Durchläuft man die Relevanz-Intervalle, so fügt jedes Intervall genau eine arithmetische Folge ($\pm$3·2·p; $p^2$) zu dem ‚Zopf' der Vielfachen-Folgen hinzu. Also gilt:
**Satz 8.7.** *p generiert die arithmetische Folge ($\pm$3·2·p; $p^2$) $\subset$ H der Vielfachen von p.*
**Satz 8.8.** *Alle Glieder von ($H_+$; $p_i^2$) oder ($H_-$; $p_i^2$), die in einer arithmetischen Folge ($\pm$3·2·$p_i$; $p_i^2$) vorkommen, sind Vielfache. Alle anderen sind Primzahlen.*

Die Glieder der Folgen ($\pm$3·2·p; $p^2$) sind $H_+$ u. $H_-$ zugeordnet. Suffix + an p sagt, p hat die Gestalt 3·2·k+1; analog Suffix -. $p_+ \cdot H_+$ bedeutet, jedes Glied von $H_+$ ist multipliziert mit $p_+$. Die Primzahlen agieren dann wie Koeffizienten. Dann gilt: (3·2·$p_+$; $p_+^2$) = $p_+ \cdot H_+$ u. (3·2·$p_-$; $p_-^2$) = $p_- \cdot H_-$.
Die Vielfachen von p, die zu $H_+$ gehören, sind Glieder von (+3·2·p; $p^2$). Denn: 1.) $p^2 \in H_+$ nach Satz 6.2. Die Glieder von (+3·2·p; $p^2$) haben die Gestalt $p^2$ + 3·2·p·k, werden also von p geteilt. 3.) (+3·2·p; $p^2$) $\in H_+$, weil $p^2 \in H_+$ u. andere Glieder haben zu $p^2$ einen Abstand von Vielfachen von 3·2. 4.) Die konstante Differenz aufeinander folgender Glieder ist 3·2·p, weil die Rest-Systeme (mod p) in $H_+$ vollständig sind. 5.) Es gibt keine anderen Glieder von (+3·2·p; $p^2$) in $H_+$, weil sie positive Reste (mod p) lassen würden. Die Glieder von ($\pm$p; 0), sofern sie zu H gehören, sind El. des Durchschnitts ($\pm$p; 0) $\cap$ H. Die El. dieses Durchschnitts sind Glieder der arithmetischen Folge p·H. H = $H_+ \cup H_-$. Also gilt:
**Satz 8.9.** *Die Vielfachen von p, sofern sie zu ($\pm$3·2·p; $p^2$) gehören, erscheinen abwechselnd in den Ästen von H. $H_+$ enthält {-5·$H_-$, +7·$H_+$, -11·$H_-$, +13·$H_+$, -17·$H_-$, ..}. $H_-$ enthält {-5·$H_+$, +7·$H_-$, -11·$H_+$, +13·$H_-$, -17·$H_+$, ..}.*

## § 9. Primzahlen in Relevanz-Intervallen.

Nächstes Ziel dieser Arbeit ist der Beweis der Unendlichkeit der Menge der Primzahlzwillinge, § 10. Dazu ist genauere Kenntnis der Primzahl-Situation in den Relevanz-Intervallen nützlich.

Sind a und b Primzahlzwillinge, d.h. a, b $\in$ P und etwa a > b und dann a – b = 2, so ist a $\in H_+$ und b $\in H_-$. Wir sagen: a und b sind Glieder von ‚parallelen' Relevanz-Intervallen ($H_+$; $p^2$) und ($H_-$; $p^2$); also gibt es dann immer ein k mit a = 3·2·k + 1 $\in$ ($H_+$; $p^2$) und b = 3·2·k – 1 $\in$ ($H_-$; $p^2$), mit gleichem k für beide. Zum Unendlichkeits-Beweis der Menge der Primzahlzwillinge ist u. a. die Frage der Existenz von Primzahlen in ‚parallelen' Relevanz-Intervallen, ($H_+$; $p^2$) und ($H_-$; $p^2$), zu diskutieren. Zum Beweis ist nicht erforderlich vorbereitend zu beweisen, dass a l l e ‚parallelen' Relevanz-Intervalle Primzahlen enthalten. Es reicht - etwas schwächer, aber leichter beweisbar - zu zeigen, dass es $\infty$ v i e l e ‚parallele' Relevanz-Intervalle gibt, die Primzahlen enthalten.

Zur Einstimmung sei auf Sätze der Literatur hingewiesen. Erstens: [3, Seite 27]: „Wenn auch die Abstände der Primzahlen sehr unregelmäßig sind, so kann man doch Aussagen über die Häufigkeit der Primzahlen im großen ganzen machen. Die Summe $\Sigma(1/p)$ der Reziproken aller Primzahlen divergiert, was EULER zum Nachweis der Existenz unendlich vieler Primzahlen verwandte. Dagegen konvergiert die Reihe $\Sigma(1/n^2)$, sogar die Reihe $\Sigma(1/n^{1+\in})$, $\in>0$, wenn über alle natürlichen n>0 summiert wird. Die Primzahlen liegen also dichter als die Quadratzahlen.“ Diese Aussage läßt sich allerdings nicht auf einzelne Relevanz-Intervalle herunter brechen. – Zweitens: „Ein prägnantes Ergebnis ist der von GAUß vermutete und von HADAMARD und DE LA VALLEE-POUSSIN im Jahre 1896 bewiesene Primzahlsatz: Das Verhältnis der Anzahl $\pi(n)$ der Primzahlen bis n und der Funktion n:log n strebt mit wachsendem n gegen1.“ Aber ein Ausdruck $\pi(n+1) - \pi(n)$ gibt, als Gegenstand einer Abschätzung, nur scheinber eine belastbare Größe an. - Drittens: „Ein Ergebnis über die Verteilung der Primzahlen auf gewisse Klassen stammt von DIRICHLET, SATZ III. Er hat den berühmten Satz über die arithmetische Progression bewiesen: in jeder arthmetischen Progression a, a+m, a+2m, a+3m, … , in der (a,m)=1 gilt, gibt es unendlich viele Primzahlen.“ „Der DIRICHLETsche Satz liefert sogar eine gewisse Gleichverteilung der Primzahlen auf die verschiedenen Progressionen $a+n\cdot m$, n = 1, 2, … und (a,m)=1.“

**Satz 9.1.** Es gibt unendlich viele Relevanz-Intervalle der Gestalt $(H_+; p^2)$ und es gibt unendlich viele Relevanz-Intervalle der Gestalt $(H_-; p^2)$.
Beweis: Denn P ist unendlich.

Aus den zitierten Sätzen können wir schließen:
**Satz 9.2.** Es gibt unendlich viele Relevanz-Intervalle in $H_+$, die Primzahlen enthlaten, und es gibt unendlich viele Relevanz-Intervalle in $H_-$, die Primzahlen enthalten.
Beweis: Denn es gibt jeweils in den Ästen unendlich viele Primzahlen – und jedes Relevanz-Intervall ist endlich.

Wir können einen wichtigen Satz hinzufügen: Satz 8.9.: Das Alternieren der Zugehörigkeit der Glieder zu $H_+$ u. $H_-$:
**Satz 9.3.** Es gibt in $(H_+; p^2)$ etwa gleich viele Vielfache wie in $(H_-; p^2)$.
Beweis: Denn das Vorkommen der aufeinander folgenden Vielfachen, die zu den Intervallen gehören, alterniert nach Satz 8.9.
**Satz 9.4.** Wenn es in $(H_+; p^2)$ Primzahlen gibt, dann auch in $(H_-; p^2)$ - und umgekehrt - und zwar etwa gleich viele.
Beweis: Folgt aus Satz 8.9.
**Satz 9.5.** Es gibt unendlich viele ‚parallele‘ Relevanz-Intervalle, die Primzahlen enthalten.
Beweis: Folgt aus Satz 9.4. - Es handelt sich dann zwangsläufig um Primzahlen beider Gestalten.

| p | q | L | $(H_+; p^2)$ | | $(H_-; q^2)$ | $(H_+; p^2)$ | | $(H_-; p^2)$ | d = q - p |
|---|---|---|---|---|---|---|---|---|---|
| | | | # Primz. | ~ | #Primz. | # Vielf. | ~ | # Vielf. | |
| 5 | 7 | 4 | 3 | | 3 | 1 | | 1 | 2 |
| 11 | 13 | 8 | 5 | | 3 | 3 | | 5 | 2 |
| 17 | 19 | 12 | 5 | | 5 | 7 | | 7 | 2 |
| 29 | 31 | 20 | 7 | | 10 | 13 | | 10 | 2 |
| 41 | 43 | 28 | 12 | | 7 | 16 | | 21 | 2 |
| 59 | 61 | 40 | 18 | | 13 | 22 | | 27 | 2 |
| 71 | 73 | 48 | 14 | | 17 | 34 | | 31 | 2 |
| 101 | 103 | 68 | 18 | | 23 | 50 | | 45 | 2 |
| 107 | 109 | 72 | 21 | | 22 | 51 | | 50 | 2 |
| 7 | 11 | 12 | 7 | | 9 | 5 | | 3 | 4 |
| 13 | 17 | 20 | 10 | | 13 | 10 | | 7 | 4 |
| 19 | 23 | 28 | 13 | | 15 | 15 | | 13 | 4 |
| 37 | 41 | 52 | 22 | | 23 | 30 | | 29 | 4 |
| 67 | 71 | 92 | 34 | | 31 | 58 | | 61 | 4 |
| 79 | 83 | 108 | 39 | | 36 | 69 | | 72 | 4 |
| 97 | 101 | 132 | 42 | | 47 | 90 | | 85 | 4 |
| 103 | 107 | 140 | 46 | | 41 | 94 | | 99 | 4 |
| 109 | 113 | 148 | 49 | | 51 | 99 | | 97 | 4 |
| 23 | 29 | 52 | 24 | | 22 | 28 | | 30 | 6 |
| 31 | 37 | 68 | 28 | | 28 | 40 | | 40 | 6 |
| 47 | 53 | 100 | 39 | | 42 | 61 | | 58 | 6 |
| 53 | 59 | 112 | 38 | | 40 | 74 | | 72 | 6 |
| 61 | 67 | 128 | 42 | | 49 | 86 | | 79 | 6 |
| 73 | 79 | 152 | 53 | | 53 | 99 | | 99 | 6 |
| 83 | 89 | 172 | 52 | | 61 | 120 | | 111 | 6 |
| … | | | | | | | | | |

Beispiele, d. h. wahre Zahlen, zu den obigen Sätzen, anfängliche Intervalle.
Sinnvoll zu vergleichen sind Intervalle mit gleichem d.

Tabelle 9.1.

## § 10. Unendlichkeit der Menge der Primzahlzwillinge.

Behauptung. Es gibt ∞ viele Primzahlzwillinge.
Zum Beweis die Aussagen A. bis H:

A. Der Geometrische Ort der Menge der Primzahlzwillinge im Primzahl-Gitter.
Wir arbeiten im Primzahl-Gitter H × H. Aus den Koordinaten auf den Halb-Achsen werden die Punkte der Folge ($H_-$, $H_+$) = (1, 1), (-5, 7), (-11, 13), (-17, 19), (-23, 25), (-29, 31), ... gebildet. Das sind im Primzahl-Gitter die Punkte der Winkelhalbierenden des 2. Quadranten. In jedem dieser Punkte ist die Summe der Koordinaten gleich +2. Siehe Tab. 10.1.:

2. Quadrant $H_+$ 1. Quadrant
↑
…
(-29,31) 31
(-23,25) 25
(-17,19) 19
(-11, 13) 13
(-5,7) 7
$H_-$ ← …-47 -41 -35 - 29 - 23 -17 - 11 - 5 (1, 1) 7 13 19 25 31 37 43 49 55 … → $H_+$
-5
-11
-17 -
-23
-29
…
↓
3. Quadrant $H_-$ 4. Quadrant

Primzahl-Gitter H × H
Winkelhalbierende im 2. Quadranten von H × H als geometrischer Ort der Zwillinge, insbes. der Primzahlzwillinge
Tabelle 10.1

Wir multiplizieren die negativen Koordinaten mit (-1). Dann sind alle Koordinaten positiv u. wir bewegen uns in der üblichen Darstellung, in der die Primzahlen positiv sind, s. Tab. 10.1.:

(Satz 4.1.:) Im Primzahl-Gitter H × H ist die Winkelhalbierende des 2. Quadranten der geometrische Ort der Menge der Zwillinge, insbes. der Menge der Primzahlzwillinge. Denn sie enthält alle Zwillinge, insbes. alle Primzahlzwillinge. Die Punkte auf $H_-$ haben die Gestalt 3·2·s-1, auf $H_+$ die Gestalt 3·2·s+1. Sie sind entweder Vielfache $v_-$ oder Primzahlen $p_-$, oder $v_+$ oder $p_+$. Die Punkte von H × H sind dann (in kurzer Bezeichnung) entweder $v_-v_+$ oder $v_-p_+$ oder $p_-v_+$ oder $p_-p_+$-Punkte; andere Kombinationen gibt es nicht.
Die Halbachsen seien wie oben in Relevanz-Intervalle ($H_-$; $p^2$) u. ($H_+$; $p^2$) eingeteilt. Beispiel: Das erste Relevanz-Intervall auf der Winkelhalbierenden ist {($H_-$; $1^2$), ($H_+$; $1^2$)} = {(-1, 1), (5, 7), (11, 13), (17, 19)}. Es beginnt mit dem Glied ‚untypischen' (-1, 1). Es enthält drei Primzahlzwillinge. Der Primzahlzwilling (3, 5) erscheint nicht, da 3 nicht zu H gehört.
B. Das Übergewicht der Vielfachen als Voraussetzung für den Beweis. P hat zwei zueinander fremde Abschnitte. Im ersten Abschnitt, p < 23, haben in den Relevanz-Intervallen die Primzahlen das Übergewicht - im zweiten, p ≥ 23, die Vielfachen. Im 1. Abschnitt gibt es nach dem Schubfächer-Prinzip (Satz II) Primzahlzwillinge. Im 2. Abschnitt gibt es - ebenfalls nach dem Schubfächer-Prinzip - Vielfachen-Zwillinge. Die Existenz von Primzahlzwillingen ist im 2. Abschnitt nicht trivial u. ist zu beweisen. Das Übergewicht der Vielfachen ist essentiell für den Beweis. Wir betrachten nur Relevanz-Intervalle im 2. Abschnitt, d.h. mit Übergewicht der Vielfachen.
Vorweg sei festgestellt: man hat keinen Einfluß darauf, welches Glied von $H_-$ mit welchem Glied von $H_+$ zu einem Punkt der Winkelhalbierenden zusammenkommt. Die Zuordnung wird ausschließlich von der Gestalt von H bestimmt. H wiederum erhält seine Gestalt von den Vielfachen-Folgen der relevanten Primzahlen. Ziel unserer Zuordnungs-Betrachtung ist, möglichst zu erkennen, ob u. welche Zwangsläufigkeiten bei der Bildung der Paare wirksam sind. Solche Zwangsläufigkeiten tragen bei zum Verständnis der S t r u k t u r der Primzahlverteilung.

C. Ausgangsituation. Beispiel zum Verfahren. Die folgenden Aussagen (1) bis (11) gelten für alle p mit p ≥ 23. Zur Erleichterung der Lektüre wird als Beispiel der Fall p = 41 in Kleindruck, (s. Tab. 10.2.), mitgeführt. Anzahlen werden mit # angegeben. Dieser Fall ist so ‚groß', dass er alle wissenswerten Eigenschaften aufweist – und so ‚klein', dass er noch in den Text passt. Wir betrachten also als Beispiel die Relevanz-Intervalle ($H_-$; $41^2$) u. ($H_+$; $41^2$). Im Sinne von §9 enthalten beide Intervalle Primzahlen.

(1) $(H_-; 41^2)$ = 1679, 1685, 1691, 1697, 1703, 1709, 1715, 1721, 1727, 1733, 1739, 1745, 1751, 1757, 1763, 1769, 1775, 1781, 1787, 1793, 1799, 1805, 1811, 1817, 1823, 1829, 1835, 1841.

$(H_+; 41^2)$ = 1681, 1687, 1693, 1699, 1705, 1711, 1717, 1723, 1729, 1735, 1741, 1747, 1753, 1759, 1765, 1771, 1777, 1783, 1789, 1795, 1801, 1807, 1813, 1819, 1825, 1831, 1837, 1843.

Die $v_-$Glieder u. die $p_-$Glieder liegen in $(H_-; 41^2)$. Die $v_+$Glieder u. die $p_+$Glieder liegen in $(H_+; 41^2)$.
Wir kennen die Anzahlen. Sie sind durch (1) gegeben:

(2) $\#v_-, \#p_-, \#v_+, \#p_+$.
21 17 16 12

In beiden Intervallen existieren Vielfache u. Primzahlen:

(3) $\#v_- > 0, \#p_- > 0, \#v_+ > 0, \#p_+ > 0$.
21 7 16 12

Durch die Wahl von $p \geq 23$ sind die Vielfachen in der Überzahl:

(4) $\#v_- > \#p_+$ u. $\#v_+ > \#p_-$.
21 > 12 16 > 7

Zu vergleichen sind die v-Minus-Glieder mit den p-Plus-Gliedern u. die v-Plus-Gl. mit den p-Minus-Gl.
Nach (3) u. (4) gibt es zwischen den Einzel-Posten die Überstände:

(5) $\#v_- - \#p_+ = 9 > 0$ u. $\#v_+ - \#p_- = 9 > 0$.
21 - 12 = 9 > 0 16 - 7 = 9 > 0

Die Intervall-Länge L ist in beiden Intervallen

(6) $L = (q^2 - p^2) / 3 \cdot 2$.
$(43^2 - 41^2) / 3 \cdot 2 = 28$

Die Summe der Anzahlen ergibt L:

(7) $\#v_+ + \#p_+ = \#v_- + \#p_- = L$.
16 + 12 = 21 + 7 = 28

Überstände vor der 1. Streichung. Nach (4) u. (7) sind diese Überstände der Einzel-Posten positiv u. gleich:

(8) $\#v_+ - \#p_- = \#v_- - \#p_+ > 0$.
16 - 7 = 21 - 12 = 9 > 0

Die geeigneten Summen der Einzelposten ergeben L:

(9) $\#v_- + \#p_- = \#v_+ + \#p_+ = \#v_-v_+ + \#v_-p_+ + \#p_-v_+ + \#p_-p_+ = L$.
21 + 7 = 16 + 12 = 12 + 4 + 9 + 3 = 28 = L

Da die Punkte auf der Winkelhalbierenden existieren u. da (4) gilt, sind die Anzahlen der Paare in den Paar-Typen $v_-v_+$, $p_-v_+$, $v_-p_+$ u. $p_-p_+$ nicht alle gleich 0. Andere Paar-Typen gibt es nicht:

(10) $\#v_-v_+ > 0$ oder $\#p_-v_+ > 0$ oder $\#v_-p_+ > 0$ oder $\#p_-p_+ > 0$. Ob $\#p_-p_+ > 0$ ist, ist Gegenstand der Betrachtung.
12 > 0 oder 9 > 0 oder 4 > 0 oder 3 > 0

Nach (3) (u. nach dem Primzahlsatz) gibt es $v_-v_+$Glieder:

(11) $\#v_-v_+ > 0$.
12 > 0

D. <u>Verkleinerung der ‚Betrachtungsmasse'</u> durch Abbau von Überständen.

Bemerkung: Die Anordnung der Primzahlen im Relevanz-Intervall scheint sehr verwickelt zu sein. Dennoch ist eine Struktur vorhanden. Es gilt, sie zu erkennen u. auszuwerten.

Überstand bedeutet: Es sind mehr 'paarungswillige' als 'paarungsbereite' Glieder vorhanden. Dann kann nicht jedes 'paarungswillige' Glied bei der Paarung ein 'paarungsbereites' Glied abbekommen. Wir sind bestrebt, die Überstände abzubauen. Das Abbauen von Überständen kann man so weit - aber auch nur so weit - treiben, wie Überstände vorhanden sind. Klar ist, dass in jedem Paar $(H_-; p^2)$ u. $(H_+; p^2)$ die Größen $\#v_-$, $\#p_-$, $\#v_+$ u. $\#p_+$ existieren. Dann ist auch klar, dass sie ein Minimum besitzen, u. dass das Minimum unter den kleineren zu suchen ist, also $\#p_-$ oder $\#p_+$. Im Normalfall sind die beiden verschieden. Wir erwarten also, dass der Abbau der Überstände bis $\min(\#p_-, \#p_+)$ geht (im Beispiel ist $\min(\#p_-, \#p_+) = \#p_- = 7$, s.Tabelle 10.2.).

Die Überstände führen zu zwangsläufigen Paarungen. Sie greifen - nach Vollzug der Zwangsläufigkeit - in das Paarungsgeschehen nicht mehr ein. Deswegen sind sie im Sinne der Betrachtung nicht aussagekräftig. Sie werden 'unspezifisch' (d.h. einerlei welche von ihrem Typ) gestrichen. Dadurch verkleinert sich die Betrachtungs-Masse.

E. Die <u>erste Streichung</u>: Überstände von vGliedern gegenüber pGliedern: v > p. Nach Voraussetzung ist $\#v_- > \#p_+$ u. $\#v_+ > \#p_-$.

In der 1. Streichung werden diejenigen Überstände behandelt, die $v_-$Glieder gegenüber $p_+$Gliedern u. $v_+$Glieder gegenüber $p_-$Gliedern aufweisen. Der schließlich angefügte Stern * bedeutet: 'Wert nach der 1. Streichung'.

Die Überstände $\#v_+ - \#p_- = \#v_- - \#p_+ > 0$ in (5) geben Anlaß zu der Aussage: die $\#v_+ - \#p_-$ (16 - 7 = 9) überständigen $v_+$Glieder haben bei der Paarbildung keine andere Möglichkeit als $v_-$Glieder zu wählen. Umgekehrt müssen die $\#v_-$

- $\#p_+$ (21 - 12 = 9) überständigen $v_-$Glieder bei der Paarbildung $v_+$Glieder wählen. Die Überstände dieser Art sind immer trivialerweise einander gleich:

(12) $\#v_- - \#p_+ = \#v_+ - \#p_- = 9 > 0$.

21 - 12 = 16 - 7 = 9 > 0

Wir erkennen (im Beispiel) 9 $v_-v_+$Glieder. D. h.: von den vorhandenen 12 $v_-v_+$Gliedern dürfen wir 9 streichen. Es ist einerlei welche. Es haben also mindestens 9 überständige $v_-v_+$Glieder ihren Paarbildungsauftrag erfüllt u. werden gestrichen. S. Tabelle 10.2.

Nach dieser 1. Streichung sind noch vorhanden:

verbleibende $\#v_+{}^*v_+$Glieder = ursprüngliche $\#v_+$ - überstehende $\#v_+$,

7 = 16 - 9

verbleibende $\#v_-{}^*v$-Glieder = ursprüngliche $\#v_-$ - überstehende $\#v_-$ .

12 = 21 - 9

Außerdem verbleiben die 12 $p_+$Glieder, da von der 1. Streichung unberührt. Ebenso verbleiben die 7 $p_-$Glieder, da ebenfalls von der 1. Streichung unberührt. Außerdem: es gibt nur noch urspüngliche $\#v_-v_+$ - gestrichene $v_-v_+$Glieder = $\#v_-v_+{}^* = \#v_-v_+ - (\#v_+ - \#p_-) = \#v_-v_+ - (\#v_- - \#p_+)$.

12 - 9 = 3 = 12 - (16 - 7) = 12 - (21 - 12)

D.h. immer:

(13) $\#v\text{-}^* = \#p_+$ u. $\#v_+{}^* = \#p_-$

= 12 = 7 .

F. Die <u>2. Streichung.</u> Überstände: $v_-$Glieder gegenüber $v_+$Gliedern, $p_+$Glieder gegenüber $p_-$Gliedern. Aussage (13) zeigt, dass bei $\#v_-{}^*$ u. $\#p_+$ noch Handlungsbedarf besteht, während $\#v_+{}^*$ u. $\#p_-$ das angekündigte Streichminimum erreicht haben.

In der 2. Streichung werden diejenigen Überstände behandelt, die $v_-$Glieder noch gegenüber $v_+$Gliedern aufweisen, sowie diejenigen Überstände, die $p_+$Glieder noch gegenüber $p_-$Gliedern aufweisen. Der angefügte Doppel-Stern ** bedeutet: 'Wert nach der 2. Streichung'. In diesem Sinne gibt es weitere Überstände von Einzelposten:

(14) $\#v_-{}^* > \#v_+{}^*$ u. $\#p_+ > \#p_-$,

Im Beispiel ist:

12 $p_+$Glieder - 7 $p_+$Glieder ergeben 5 überständige $p_+$Glieder. Sie finden keine $p_-$Glieder.

12 $v_-$Glieder - 7 $v_-$Glieder ergeben 5 überständige $v_-$Glieder. Sie finden keine $v_+$Glieder.

Diese 5 überständigen $p_+$Glieder u. $v_-$Glieder finden keine $p_-$Glieder u. $v_+$Glieder. Sie bilden 5 $p_+v_-$Glieder. Auch diese Glieder haben ihren Paarbildungsauftrag erfüllt. Sie werden gestrichen. Siehe Tabelle 10.2. Nach dieser 2. Streichung gibt es noch

$\#v_-{}^{**} = 7$, $\#p_+{}^{**} = 7$

u.

$\#v_+{}^* = 7$, verblieben nach der 1. Streichung u. wurden nicht weiter abgebaut,

$\#p_- = 7$, als Minimum von beiden Streichungen unberührt.

D.h.:

(15) $\#v_-{}^{**} = \#p_- = \#v_+{}^* = \#p_+{}^{**} = 7$. Nach der 2. Streichung gibt es keine Überstände von Einzel-Posten. Und die vier verbleibenden Einzel-Posten sind einander gleich. Das 'Streich-Minimum' ist erreicht.

Außerdem gibt es nur noch

9 - 5 $p_-v_+$Glieder = 4 $p_-v_+$Glieder: $\#p_-v_+{}^* = 4$. Und ebenso ist $\#v_-p_+{}^* = 4$.

Nach beiden Streichungen sind noch vorhanden:

(16) $\#v_-{}^{**} = \#p_- = \#v_+{}^* = \#p_+{}^{**}$ (= 7 im Beispiel), (= Wiederholung von (15).

(17) $\#v_-v_+{}^{**} = \#p_-p_+$ (= 3 im Beispiel) u. $\#v_-p_+{}^{**} = \#p_-v_+{}^{**}$ (= 4 im Beispiel).

D.h.:

**Satz 10.1.** Wenn $\#v_-v_+{}^{**} > 0$ ist, so gilt auch $\#p_-p_+ > 0$ - und umgekehrt.

Das Junctim in Satz 10.1. ist ein Kriterium für das gesuchte Nicht-Verschwinden von $\#p_-p_+$.

G. $\#v_-v_+{}^{**} > 0$. Häufigkeiten.

Auch die Häufigkeiten besitzen eine extreme Besonderheit.

Wir knüpfen an Aussage (17) oder Satz 10.2. an. Damit #pp = 0 ist, muß #pv = #vp (im Beispiel = 7) sein. D. h. alle restlichen p müssen ein restliches v finden. Diesen Fall gibt es nur <u>1 Mal</u>. Die Anzahl dieses Falles sei a und es ist

a = 1. Aber für die anderen Fälle gibt es eine hohe Zahl von Kombinationen. Die Anzahl Kombiationen dieser Fälle sei b.

(18) Der Fall 'alle restlichen p finden bei der Paarbildung ein restliches v' ist einzig. Er bleibt einzig auch bei wachsendem p u. gilt für alle p: a = 1 für alle p ≥ 23.

(19) Bei wachsendem p ist b eine schnell wachsende große Zahl. b wächst über alle Grenzen, während a in allen Intervallen = 1 bleibt: b >> a und a = 1 für alle p.

Die Häufigkeiten zeigen, dass im Normalfall d. h. in unendlich vielen Fällen, $\#v_-v_+ > 0$ ist. Also ist im Normalfall $\#p_-p_+ > 0$. q. e. d. Also gilt:

(20) Es gibt ∞ viele Intervalle in denen #vv u. #pp nicht verschwinden.

L wächst mit p nicht streng monoton. Wegen des Abnehmens der Primzahl-Dichte nach dem Primzahlsatz ist aber der Mangel an Monotonie nicht von Bedeutung für unsere Betrachtung. Die Aussagen (1) bis (20) und Satz 10.2. gelten - mutatis mutandis - für alle p ≥ 23. Damit ist gezeigt:

**Satz 10.2.** *Es gibt ∞ viele Primzahlzwillinge.*

Der Beweis hat ein Nebenergebnis zur Struktur erbracht:

**Satz 10.3.** *Ist in parallelen Relevanz-Intervallen #v > #p, s gibt es dort vv-, vp-, pv- und pp-Paare.*

Beispiel p = 41

| 1. Schritt originäre Werte feststellen Paar-Typ bestimmen | | | | 2. Schritt Paare sortieren nach Paar-Typ | | | | 3. Schritt 1. Überstand zu vv | | | | 4. Schritt 2. Überstand zu pv oder vp | | | | 5. Schritt verbleibende sichten | | | |
|---|---|---|---|---|---|---|---|---|---|---|---|---|---|---|---|---|---|---|---|
| Typ | Nr | $(H_-, 41)$ | $(H_+, 41)$ | Typ | Nr | $(H_-, 41)$ | $(H_+, 41)$ | Typ | Nr | $(H_-, 41)$ | $(H_+, 41)$ | Typ | Nr | $(H_-, 41)$ | $(H_+, 41)$ | Typ | Nr | $(H_-, 41)$ | $(H_+, 41)$ |
| vv | 1 | 1679 | 1681 | vv | 9 | 1727 | 1729 | ~~vv~~ | ~~9~~ | ~~1727~~ | ~~1729~~ | | | | | | | | |
| vv | 2 | 1685 | 1687 | vv | 7 | 1715 | 1717 | ~~vv~~ | ~~7~~ | ~~1715~~ | ~~1717~~ | | | | | | | | |
| vp | 3 | 1691 | 1693 | vv | 5 | 1703 | 1705 | ~~vv~~ | ~~5~~ | ~~1703~~ | ~~1705~~ | | | | | | | | |
| pp | 4 | 1697 | 1699 | vv | 28 | 1841 | 1843 | ~~vv~~ | ~~28~~ | ~~1841~~ | ~~1843~~ | | | | | | | | |
| vv | 5 | 1703 | 1705 | vv | 27 | 1835 | 1837 | ~~vv~~ | ~~27~~ | ~~1835~~ | ~~1837~~ | | | | | | | | |
| pv | 6 | 1709 | 1711 | vv | 24 | 1817 | 1819 | ~~vv~~ | ~~24~~ | ~~1817~~ | ~~1819~~ | | | | | | | | |
| vv | 7 | 1715 | 1717 | vv | 22 | 1805 | 1807 | ~~vv~~ | ~~22~~ | ~~1805~~ | ~~1807~~ | | | | | | | | |
| pp | 8 | 1721 | 1723 | vv | 20 | 1793 | 1795 | ~~vv~~ | ~~20~~ | ~~1793~~ | ~~1795~~ | | | | | | | | |
| vv | 9 | 1727 | 1729 | vv | 2 | 1685 | 1687 | ~~vv~~ | ~~2~~ | ~~1685~~ | ~~1687~~ | | | | | | | | |
| pv | 10 | 1733 | 1735 | vv | 16 | 1769 | 1771 | vv | 16 | 1769 | 1771 | vv | 16 | 1769 | 1771 | vv | 16 | 1769 | 1771 |
| vp | 11 | 1739 | 1741 | vv | 15 | 1763 | 1765 | vv | 15 | 1763 | 1765 | vv | 15 | 1763 | 1765 | vv | 15 | 1763 | 1765 |
| vp | 12 | 1745 | 1747 | vv | 1 | 1679 | 1681 | vv | 1 | 1679 | 1681 | vv | 1 | 1679 | 1681 | vv | 1 | 1679 | 1681 |
| vp | 13 | 1751 | 1753 | pv | 6 | 1709 | 1711 | pv | 6 | 1709 | 1711 | pv | 6 | 1709 | 1711 | pv | 6 | 1709 | 1711 |
| vp | 14 | 1757 | 1759 | pv | 25 | 1823 | 1825 | pv | 25 | 1823 | 1825 | pv | 25 | 1823 | 1825 | pv | 25 | 1823 | 1825 |
| vv | 15 | 1763 | 1765 | pv | 23 | 1811 | 1813 | pv | 23 | 1811 | 1813 | pv | 23 | 1811 | 1813 | pv | 23 | 1811 | 1813 |
| vv | 16 | 1769 | 1771 | pv | 10 | 1733 | 1735 | pv | 10 | 1733 | 1735 | pv | 10 | 1733 | 1735 | pv | 10 | 1733 | 1735 |
| vp | 17 | 1775 | 1777 | vp | 3 | 1691 | 1693 | vp | 3 | 1691 | 1693 | ~~vp~~ | ~~3~~ | ~~1691~~ | ~~1693~~ | | | | |
| vp | 18 | 1781 | 1783 | vp | 26 | 1829 | 1831 | vp | 26 | 1829 | 1831 | ~~vp~~ | ~~26~~ | ~~1829~~ | ~~1831~~ | | | | |
| pp | 19 | 1787 | 1789 | vp | 21 | 1799 | 1801 | vp | 21 | 1799 | 1801 | ~~vp~~ | ~~21~~ | ~~1799~~ | ~~1801~~ | | | | |
| vv | 20 | 1793 | 1795 | vp | 18 | 1781 | 1783 | vp | 18 | 1781 | 1783 | ~~vp~~ | ~~18~~ | ~~1781~~ | ~~1783~~ | | | | |
| vp | 21 | 1799 | 1801 | vp | 17 | 1775 | 1777 | vp | 17 | 1775 | 1777 | ~~vp~~ | ~~17~~ | ~~1775~~ | ~~1777~~ | | | | |
| vv | 22 | 1805 | 1807 | vp | 14 | 1757 | 1759 | vp | 14 | 1757 | 1759 | vp | 14 | 1757 | 1759 | vp | 14 | 1757 | 1759 |
| pv | 23 | 1811 | 1813 | vp | 13 | 1751 | 1753 | vp | 13 | 1751 | 1753 | vp | 13 | 1751 | 1753 | vp | 13 | 1751 | 1753 |
| vv | 24 | 1817 | 1819 | vp | 12 | 1745 | 1747 | vp | 12 | 1745 | 1747 | vp | 12 | 1745 | 1747 | vp | 12 | 1745 | 1747 |
| pv | 25 | 1823 | 1825 | vp | 11 | 1739 | 1741 | vp | 11 | 1739 | 1741 | vp | 11 | 1739 | 1741 | vp | 11 | 1739 | 1741 |
| vp | 26 | 1829 | 1831 | pp | 8 | 1721 | 1723 | pp | 8 | 1721 | 1723 | pp | 8 | 1721 | 1723 | pp | 8 | 1721 | 1723 |
| vv | 27 | 1835 | 1837 | pp | 4 | 1697 | 1699 | pp | 4 | 1697 | 1699 | pp | 4 | 1697 | 1699 | pp | 4 | 1697 | 1699 |
| vv | 28 | 1841 | 1843 | pp | 19 | 1787 | 1789 | pp | 19 | 1787 | 1789 | pp | 19 | 1787 | 1789 | pp | 19 | 1787 | 1789 |
| #v- = 21 | | | | Sortierung Paar-Typ | | | | 1. Streichung: | | | | 2. Streichung: | | | | $\#v_- = 7$ | | | |
| #p- = 7 | | | | | | | | $\#v_- - \#p_+ = \#v_+ - \#p- =$ | | | | $\#v_- - \#v_+ = 21 - 16 = 5$ | | | | #p- = 7 | | | |
| $\#v_+ = 16$ | | | | | | | | 21 -12 = 16 - 7 = 9 | | | | $\#p_+ - \#p_- = 12 - 7 = 5$ | | | | $\#v_+ = 7$ | | | |
| $\#p_+ = 12$ | | | | | | | | 9 $v_-v_+$ sind zu streichren | | | | 5 v-$p_+$ sind zu streichen | | | | $\#p_+ = 7$ | | | |
| | | | | $\#v_-v_+ = 14$ | | | | | | | | | | | | $\#v_-v_+ = 3$ | | | |
| | | | | $\#p_-p_+ = 7$ | | | | | | | | | | | | $\#p_-p_+ = 3$ | | | |
| | | | | $\#p_-v_+ = 16$ | | | | | | | | | | | | $\#p_-v_+ = 4$ | | | |
| | | | | $\#v_-p_+ = 15$ | | | | | | | | | | | | $\#v_-p_+ = 4$ | | | |

$\#v_- + \#p_- = \#v_+ + p_+ = (43^2 - 41^2) / 2 \cdot 3 = 28$

$\#v- -\# p_+ = \#v_+ - \#p- = 21 - 12 = 16 - 7 = 9$

$\#v_- - \#v+ = \#p_+ - \#p_- = 5$

$\#v_-v_+ - \#p_-p_+ = 14 - 7 = 7$

$\#v_- - \# p_+ = \#v_+ - \#p_- = \#v_-v_+ - \#p_-p_+ = 21 - 12 = 16 - 7 = 14 - 7 = 7$

Übersicht über die unspezifische Zuordnungspraxis im Falle p = 41

Tabelle 10.2.

**Satz 10.4.** *Die Beweis-Schritte 1 bis 20 funktionieren bei jeder Primzahl p ≥ 23. Dadurch ist eine ‚Struktur' der Primzahlverteilung beschrieben.*

## § 11. Bestätigung der GOLDBACH-Vermutung

Die GOLDBACH-Vermutung, oben (2GH) genannt, besagt, dass sich jede gerade Zahl > 2 als Summe zweier Primzahlen darstellen lasse. Die Lösungsmenge wurde in § 3 angegeben. Sie ist Teilmenge eines Systems von konzentrischen Rhomben. Beispiel Siehe z.B. Tabelle 11.1. Dabei ist ein Rhombus jeweils ‚zuständig' für drei gerade Zahlen – im Beispiel 16, 18, 20, 18. Das System der Lösungs-Rhomben enthält a l l e Lösungen der (2GH).

2. Quadrant: $E_0$ 1. Quadrant: $E_{+2}$

| | | | | | | | | | | | | | | | | | | | | | | |
|---|---|---|---|---|---|---|---|---|---|---|---|---|---|---|---|---|---|---|---|---|---|---|
| | | | | | | | | | | | ↑ | | | | | | | | | | | |
| | 120 | . | 108 | 102 | . | 90 | 84 | 78 | 72 | 66 | 61 | 68 | 74 | 80 | . | 92 | 98 | 104 | . | . | 122 | |
| | . | . | . | . | . | . | . | . | . | . | 55 | . | . | . | . | . | . | . | . | . | . | |
| | . | . | . | . | . | . | . | . | . | . | 49 | . | . | . | . | . | . | . | . | . | . . | |
| | 102. | | 90 | 84 | . | 72 | 66 | 60 | 54 | 48 | 43 | 50 | 56 | 62 | . | 74 | 80 | 86 | . | . | 104 | |
| | 96 | . | 84 | 78 | . | 66 | 60 | 54 | 48 | 42 | 37 | 44 | 50 | 56 | . | 68 | 74 | 80 | . | . | 98 | |
| | 90 | . | 78 | 72 | . | 60 | 54 | 48 | 42 | 36 | 31 | 38 | 44 | 50 | . | 62 | 68 | 74 | . | . | 92 | |
| | . | . | . | . | . | . | . | . | . | . | 25 | . | . | . | . | . | . | . | . | . | . | |
| | 78 | . | 66 | 60 | . | 48 | 42 | 36 | 30 | 24 | 19 | 26 | 32 | 38 | . | 50 | 56 | 62 | . | . | 80 | |
| | 72 | . | 60 | 54 | . | 42 | 36 | 30 | 24 | 18 | 13 | 20 | 26 | 32 | . | 44 | 50 | 56 | . | . | 74 | |
| | 66 | . | 54 | 48 | . | 36 | 30 | 24 | 18 | 12 | 7 | 14 | 20 | 25 | . | 37 | 43 | 49 | . | . | 67 | |
| -1·($H_-$) ← | 59 | 55 | 47 | 41 | 35 | 29 | 23 | 17 | 11 | 5 | 1 | 7 | 13 | 19 | 25 | 31 | 37 | 43 | 49 | 55 | 61 | → $H_+$ |
| | 64 | . | 52 | 46 | . | 34 | 28 | 22 | 16 | 10 | 5 | 12 | 18 | 24 | . | 36 | 42 | 48 | . | . | 66 | |
| | 70 | . | 58 | 52 | . | 40 | 34 | 28 | 22 | 16 | 11 | 18 | 24 | 30 | . | 42 | 48 | 54 | . | . | 72 | |
| | 76 | . | 64 | 58 | . | 46 | 40 | 34 | 28 | 22 | 17 | 24 | 30 | 36 | . | 48 | 54 | 60 | . | . | 78 | |
| | 82 | . | 70 | 64 | . | 52 | 46 | 40 | 34 | 28 | 23 | 30 | 36 | 42 | . | 54 | 60 | 66 | . | . | 84 | |
| | 88 | . | 76 | 70 | . | 58 | 52 | 46 | 40 | 34 | 29 | 36 | 42 | 48 | . | 60 | 66 | 72 | . | . | 90 | |
| | . | . | . | . | . | . | . | . | . | . | 35 | | . | . | . | . | . | . | . | . | . | |
| | 100 | . | 88 | 82 | . | 70 | 64 | 58 | 52 | 46 | 41 | 48 | 54 | 60 | . | 72 | 78 | 84 | . | . | 102 | |
| | 106 | . | 94 | 88 | . | 76 | 70 | 64 | 58 | 52 | 47 | 54 | 60 | 66 | . | 78 | 84 | 90 | . | . | 108 | |
| | 112 | . | 100 | 94 | . | 82 | 76 | 70 | 64 | 58 | 53 | 60 | 66 | 72 | . | 84 | 90 | 96 | . | . | 114 | |
| | 118 | . | 106 | 100 | . | 88 | 82 | 76 | 70 | 64 | 59 | 66 | 72 | 78 | . | 90 | 96 | 102 | . | . | 120 | |
| | | | | | | | | | | | ↓ | | | | | | | | | | | |

3. Quadrant: $E_{-2}$ 4. Quadrant: $E_0$

Beispiel aus der Lösungsmenge der Goldbach-Vermutung
Beispiel markiert: 16(=5+11), 18(=7+11), 20(=7+13), 18 bilden einen Lösungs-Rhombus
Tabelle 11.1.

**Definition 11. 1.** *Ein Achsen-Abschnitt $H_+(A)^{+1}$ ist ein Anfangsstück $\{h_1, h_2, h_3, \ldots, h_A\}$ von $H_+$ der Länge A, das bei Exponent +1 steigend durchlaufen wird, bei Exponent -1 fallend. Analog in den drei andereni Quadranten.*

**Definition 11.2.** *Eine Rhombusseite ($H_+(A)^{+1}$, $H_+(A)^{-1}$) wird aufgespannt von den Achsen-Abschnitten $H_+(A)^{+1}$und $H_+(A)^{-1}$). Ihre Punkte sind die Summen der aufspannenden Koordinaten. Analog bei den drei anderen Seiten des Rhombus.*

Im Beispiel der Tabelle 11.1 erscheint im 1. Quadranten die Rhombusseite ($H_+(4)^{+1}$, $H_+(4)^{-1}$) = ({1, 7, 13, 19}, {19, 13, 7, 1}) mit den Punkten {(1, 19), (7,13), (13, 7), (19, 1)}. Sie alle haben die Komponenten-Summe 20; zwei von ihnen sind Lösungen der (2GH). (Bemerkung: Wäre auch 1 eine Primzahl, würde das die Lösungs-Szenerie nicht stören - sondern bereichern.) - Der 1. Quadrant ist 'zuständig' für die Komponenten-Summen d.h. für die geraden Zahlen x = 2+3•2•k = 2, 8, 14, 20, 26, 32, ..., k = 0, 1, 2, ... . Der 3. Quadrant ist zuständig für x = 4+3•2•k = 4, 10, 16, 22, 28, 34, ..., k = 0, 1, 2, .... Der 2. und der 4. Quadrant sind zuständig für x = 6+3•2•k = 6, 12, 18, 24, 30, 36, ..., k = 0, 1, 2, ... . Die Symmetrie bzw. Unsymmetrie entsteht durch die Gestalt der Äste von H.

Betrachtet wird - zunächst im 1. Quadranten - eine beliebige Rhombusseite ($H_+(A)^{+1}$, $H_+(A)^{-1}$) aus der Lösungsmenge der (2GH). Die Seiten werden aufgespannt von Komponenten-Abschnitten der Länge A. Alle Größen beziehen sich auf A, so daß die Angabe ...(A) entfallen kann. Es gibt in den Komponenten-Abschnitten nur Primzahlen und Vielfache. Die Anfangsabschnitte $H_+^{+1}$ und $H_+^{-1}$, enthalten #p Primzahlen und #v Vielfache. Die beiden Komponenten-Abschnitte enthalten die gleichen Glieder, werden aber entgegesetzt durchlaufen. Zu zeigen ist, daß jede Rhombus-Seite pp-Punkte enthält.

Die Schritte werden durchnummeriert. Es gilt:
(1) 0 < #p < A, denn beide Abschnitte sind Anfangs-Abschnitte.
(2) 0 < #v < A, denn die Abschnitte enthalten nicht nur Primzahlen.
(3) A = #p + #v, denn jedes Glied ist entweder Primzahl oder Vielfaches.
In der Startroutine der Primzahlverteilung, definiert durch #p > #v, gibt es triviale Lösungen nach dem DIRICHLETschen Schubfachprinzip. Um der Trivialität zu entgehen, wird - durch Verlängerung von A - vorausgesetzt:
(4) #p < #v.
Die Glieder auf den Achsen sind entweder Vielfache oder Primzahlen, (3). Deswegen und wegen (1) sind #v und #p definiert und positiv. Bei den symmetrischen Paaren der Rhombusseite ($H_+(A)^{+1}$, $H_+(A)^{-1}$) gibt es dann nur die Kombinationen pp, pv, vp und vv. Es ist also sinnvoll, die Zählfunktionen #pp, #pv, #vp und #vv zu betrachten.
Aus (4) folgt nach dem Schubfachprinzip, daß es vv-Fälle gibt, so daß #vv nicht verschwindet. Ausserdem ist #vv kleiner als A, weil sonst #p verschwinden würde, was durch (1) ausgeschlossen ist. Es gilt sogar wegen des Schubfachprinzips: #v - #p ≤ #vv < A. Denn die wegen (4) überstehenden v-Glieder müssen vv-Glieder bilden, natürlich höchstens #vv viele. D.h.:
(5) 0 < #v - #p ≤ #vv < L.
Ein Verschwinden von #pv würde bedeuten, daß alle p sich mit p verbinden - und dann auch alle v mit v. Das kann nicht sein wegen der Gegenläufigkeit der Komponenten-Abschnitte. Dafür sorgen von Primzahl zu Primzahl unterschiedliche Länge der Reste-Systeme und der Primzahlsatz. Ebenso wenig kann #pv gleich A sein. Das gleiche gilt für #vp. In symmetrischen Quadranten sind die beiden sogar gleich:
(6) 0 < #pv < A, 0 < #vp < A, und in den symmetrischen Quadranten 1 und 3 gilt #pv = #vp.
Ebenso kann nicht #pp gleich A sein. Denn dann gäbe es in dem Komponenten-Abschnitt keine Vielfachen. Ob allerdings #pp verschwinden kann, ist die Frage von GOLDBACH, die hier untersucht wird:
(7) 0 ≤ #pp < A.
Die Zählfunktionen liegen also alle zwischen 0 und A - nur bei #pp ist das Verschwinden noch nicht ausgeschlossen und ist Gegenstand der Betrachtung. Allgemein gilt:
(8) #vv + #vp = #v und #pp + #pv = #p, denn der Ausdruck #vv + #vp enthält alle v und der Ausdruck #pp + #pv enthält all p.
Also gilt #vp = #v - #vv und #pv = #p - #pp und da #pv = #vp im symmetrischen 1. Quadranten ist, gilt #v - #vv = #p - #pp, d. h.
(9) #v - #p = #vv - #pp.
Die Betrachtung von (5) u. (9) zeigt: nur wenn #v - #p = #vv ist, kann #pp = 0 sein. Ist aber #v - #p < #vv, dann ist zwangsläufig #pp > 0.
Zur Bestätigung der (2GH) - im 1. Quadranten - ist also
(10) #v - #p < #vv
zu beweisen.
Bemerkung: Die für die (2GH) typische Situation besteht darin, daß man die Lösungen, die zu einer bestimmten geraden Zahl gehören, nur bis zur Rhombus-Seite verfolgen kann. Aber welche Lösungen das sind, ist nicht in eine geschlossene Form zu bringen. Denn auf jeder Rhombus-Seite herrschen neue, andere Primzahl-Situationen. Es gilt immer die Aussage (9) mit positiven Gliedern, so daß die Aussage (10) immer richtig ist. Allerdings ist es für den Beweis von (10) nicht erforderlich, die Lösungen genau anzugeben: es reicht zu zeigen, daß welche existieren.
Wir verfolgen die Beziehungen zwischen den beteiligten Größen weiter und bilden rechnerische Ausdrücke für #pp und für #vv:
Per Definition gilt: A = #pp + #pv + #vp + #vv. A ist die Anzahl der symmetrischen Paare der Rhombusseite ($H_+(A)^{+1}$, $H_+(A)^{-1}$). Es werden von A die Paare vom Typ pv und vom Typ vp gestrichen. Es verbleiben A' Paare: A' = A – (#pv + #vp). A' zählt die Paare vom Typ pp oder Typ vv. Von A' werden #v - #p überständige Paare vom Typ vv getrichen. Es verbleiben A'' Glieder: A'' = A' – (#v - #p). Die A'' Paare sind vom Typ pp oder Typ vv und bestehen zu gleichen Teilen aus p-Gliedern und v-Gliedern. Also enthalten die verbliebenen A'' Glieder zur einen Hälfte pp-Fälle und zur anderen Hälfte vv-Fälle, denn pv- und vp-Fälle sowie überständige vv-Fälle sind nicht mehr vorhanden. Also ist:
(11) #pp = A"/2 = [A -(#pv + #vp) - (#v - #p)] / 2.
Außerdem ist dann:

(12) #vv = A"/2 + (#v-#p) = [A -(#pv + #vp) - (#v - #p)] / 2 + [#v - #p].

**Satz 11.1.**

*#pp = [A -(#pv + #vp) - (#v - #p)] / 2.*

*#vv = [A -(#pv + #vp) - (#v - #p)] / 2 + [#v - #p].*

A ist die Länge des Komponenten-Abschnitts.

Wir betrachten die relativen Häufigkeiten. Dazu wird z.B. #vv abgeschätzt. Der Anteil v in A beträgt etwa #v / A und ist in beiden Abschnitten gleich. Also ist

(13) #vv ~ (#v / A)·(#v / A)·A.

Zu beweisen ist #v - #p < #vv:

(14) Mit (13) gilt:

(15) #v - #p < #vv ~ (#v / A)·(#v / A)·A = (#v)$^2$ / A.

Da #v > 0 darf durch #v geteilt werden:

(16) (#v - #p) / #v < #v / A. Oder, weil #v = A - #p ist:

(17) (#v - #p) / #v < (A - #p) / A. D.h.:

(18) (#v / #v) - (#p / #v) < A / A - #p / A. D.h.:

(19) - #p / #v < - #p / A. Also

(20) #p / A < #p / #v.

Das ist richtig, weil #v < A.

Also ist die Behauptung #v - #p < #vv richig und es gilt:

(21) 0 < #pp.

Damit ist die (2GH) im 1. Quadranten bestätigt.

Die Überlegung funktioniert auch - mutatis mutandis - für die Fälle #p > #v, liefert dann aber 0 < #vv. Beispiel:

A = 200, #p = 95, #v = 105,

#pp = 40, #pv + #vp = 110, #vv = 50.

#pp = [A -(#pv + #vp) - (#v - #p)] / 2 = [200 - 110 - 10] / 2 = 40.

#vv = [A -(#pv + #vp) - (#v - #p)] / 2 + [#v - #p] = 40 + 10 = 50.

Nr (H+(200)$^{+1}$, H+(200)$^{-1}$)

Originale

| Nr | | | pp | pv | vp | vv |
|---|---|---|---|---|---|---|
| 001 | **7** | **1201** | pp | | | |
| 002 | **13** | 1195 | | pv | | |
| 003 | **19** | 1189 | | pv | | |
| 004 | 25 | 1183 | | | | vv |
| 005 | **31** | 1177 | | pv | | |
| 006 | **37** | **1171** | pp | | | |
| 007 | **43** | 1165 | | pv | | |
| 008 | 49 | 1159 | | | | vv |
| 009 | 55 | **1153** | | | vp | |
| 010 | **61** | 1147 | | pv | | |
| 011 | **67** | 1141 | | pv | | |
| 012 | **73** | 1135 | | pv | | |
| 013 | **79** | **1129** | pp | | | |
| 014 | 85 | **1123** | | | vp | |
| 015 | 91 | **1117** | | | vp | |
| 016 | **97** | 1111 | | pv | | |
| 017 | **103** | 1105 | | pv | | |
| 018 | **109** | 1099 | | pv | | |
| 019 | 115 | **1093** | | | vp | |
| 020 | 121 | **1087** | | | vp | |
| 021 | **127** | 1081 | | pv | | |
| 022 | 133 | 1075 | | | | vv |
| 023 | **139** | **1069** | pp | | | |
| 024 | 145 | **1063** | | | vp | |

Nr (H+(200)$^{+1}$, H+(200)$^{-1}$)

sortiert nach Typ
pv und vp gestrichen

| Nr | | | pp | pv | vp | vv |
|---|---|---|---|---|---|---|
| 001 | **7** | **1201** | pp | | | |
| 006 | **37** | **1171** | pp | | | |
| 013 | **79** | **1129** | pp | | | |
| 023 | **139** | **1069** | pp | | | |
| 026 | **157** | **1051** | pp | | | |
| 033 | **199** | **1009** | pp | | | |
| 035 | **211** | **997** | pp | | | |
| 040 | **241** | **967** | pp | | | |
| 045 | **271** | **937** | pp | | | |
| 055 | **331** | **877** | pp | | | |
| 058 | **349** | **859** | pp | | | |
| 063 | **379** | **829** | pp | | | |
| 066 | **397** | **811** | pp | | | |
| 070 | **421** | **787** | pp | | | |
| 073 | **439** | **769** | pp | | | |
| 076 | **457** | **751** | pp | | | |
| 083 | **499** | **709** | pp | | | |
| 091 | **547** | **661** | pp | | | |
| 096 | **577** | **631** | pp | | | |
| 100 | **601** | **607** | pp | | | |
| ~~002~~ | ~~**13**~~ | ~~1195~~ | | ~~pv~~ | | |
| ~~003~~ | ~~**19**~~ | ~~1189~~ | | ~~pv~~ | | |
| ~~005~~ | ~~**31**~~ | ~~1177~~ | | ~~pv~~ | | |
| ~~007~~ | ~~**43**~~ | ~~1165~~ | | ~~pv~~ | | |

Nr (H+(200)$^{+1}$, H+(200)$^{-1}$)

sortiert nach Typ
ohne pv und vp und
Überstand vv gestrichen

| Nr | | | pp | pv | vp | vv |
|---|---|---|---|---|---|---|
| 001 | **7** | **1201** | pp | | | |
| 006 | **37** | **1171** | pp | | | |
| 013 | **79** | **1129** | pp | | | |
| 023 | **139** | **1069** | pp | | | |
| 026 | **157** | **1051** | pp | | | |
| 033 | **199** | **1009** | pp | | | |
| 035 | **211** | **997** | pp | | | |
| 040 | **241** | **967** | pp | | | |
| 045 | **271** | **937** | pp | | | |
| 055 | **331** | **877** | pp | | | |
| 058 | **349** | **859** | pp | | | |
| 063 | **378** | **829** | pp | | | |
| 066 | **397** | **811** | pp | | | |
| 070 | **421** | **787** | pp | | | |
| 073 | **439** | **769** | pp | | | |
| 076 | **457** | **751** | pp | | | |
| 083 | **499** | **709** | pp | | | |
| 091 | **547** | **661** | pp | | | |
| 096 | **577** | **631** | pp | | | |
| 100 | **601** | **607** | pp | | | |

| | | | | | | | | | | | | | | |
|---|---|---|---|---|---|---|---|---|---|---|---|---|---|---|
| 025 | **151** | 1057 | | pv | | | ~~010~~ | ~~**61**~~ | ~~1147~~ | ~~pv~~ | | | | |
| 026 | **157** | **1051** | pp | | | | ~~011~~ | ~~**67**~~ | ~~1141~~ | ~~pv~~ | | | | |
| 027 | **163** | 1045 | | pv | | | ~~012~~ | ~~**73**~~ | ~~1135~~ | ~~pv~~ | | | | |
| 028 | 169 | **1039** | | | vp | | ~~016~~ | ~~**97**~~ | ~~1111~~ | ~~pv~~ | | | | |
| 029 | 175 | **1033** | | | vp | | ~~017~~ | ~~**103**~~ | ~~1105~~ | ~~pv~~ | | | | |
| 030 | **181** | 1027 | | pv | | | ~~018~~ | ~~**109**~~ | ~~1099~~ | ~~pv~~ | | | | |
| 031 | 187 | **1021** | | | vp | | ~~021~~ | ~~**127**~~ | ~~1081~~ | ~~pv~~ | | | | |
| 032 | **193** | 1015 | | pv | | | ~~026~~ | ~~**157**~~ | ~~1057~~ | ~~pv~~ | | | | |
| 033 | **199** | **1009** | pp | | | | ~~027~~ | ~~**163**~~ | ~~1045~~ | ~~pv~~ | | | | |
| 034 | 205 | 1003 | | | | vv | ~~030~~ | ~~**181**~~ | ~~1027~~ | ~~pv~~ | | | | |
| 035 | **211** | **997** | pp | | | | ~~032~~ | ~~**193**~~ | ~~1015~~ | ~~pv~~ | | | | |
| 036 | 217 | **991** | | | vp | | ~~037~~ | ~~**223**~~ | ~~985~~ | ~~pv~~ | | | | |
| 037 | **223** | 985 | | pv | | | ~~038~~ | ~~**229**~~ | ~~979~~ | ~~pv~~ | | | | |
| 038 | **229** | 979 | | pv | | | ~~046~~ | ~~**277**~~ | ~~931~~ | ~~pv~~ | | | | |
| 039 | 235 | 973 | | | | vv | ~~047~~ | ~~**283**~~ | ~~925~~ | ~~pv~~ | | | | |
| 040 | **241** | **967** | pp | | | | ~~051~~ | ~~**307**~~ | ~~901~~ | ~~pv~~ | | | | |
| 041 | 247 | 961 | | | | vv | ~~052~~ | ~~**313**~~ | ~~895~~ | ~~pv~~ | | | | |
| 042 | 253 | 955 | | | | vv | ~~056~~ | ~~**337**~~ | ~~871~~ | ~~pv~~ | | | | |
| 043 | 259 | 949 | | | | vv | ~~061~~ | ~~**367**~~ | ~~841~~ | ~~pv~~ | | | | |
| 044 | 265 | 943 | | | | vv | ~~062~~ | ~~**373**~~ | ~~835~~ | ~~pv~~ | | | | |
| 045 | **271** | **937** | pp | | | | ~~068~~ | ~~**409**~~ | ~~799~~ | ~~pv~~ | | | | |
| 046 | **277** | 931 | | pv | | | ~~072~~ | ~~**433**~~ | ~~775~~ | ~~pv~~ | | | | |
| 047 | **283** | 925 | | pv | | | ~~077~~ | ~~**463**~~ | ~~745~~ | ~~pv~~ | | | | |
| 048 | 289 | **919** | | | vp | | ~~081~~ | ~~**487**~~ | ~~721~~ | ~~pv~~ | | | | |
| 049 | 295 | 913 | | | | vv | ~~087~~ | ~~**523**~~ | ~~685~~ | ~~pv~~ | | | | |
| 050 | 301 | **907** | | | vp | | ~~090~~ | ~~**541**~~ | ~~667~~ | ~~pv~~ | | | | |
| 051 | **307** | 901 | | pv | | | ~~095~~ | ~~**571**~~ | ~~637~~ | ~~pv~~ | | | | |
| 052 | **313** | 895 | | pv | | | ~~009~~ | ~~55~~ | ~~**1153**~~ | ~~vp~~ | | | | |
| 053 | 319 | 889 | | | | vv | ~~014~~ | ~~85~~ | ~~**1123**~~ | ~~vp~~ | | | | |
| 054 | 325 | **883** | | | vp | | ~~015~~ | ~~91~~ | ~~**1117**~~ | ~~vp~~ | | | | |
| 055 | **331** | **877** | pp | | | | ~~019~~ | ~~115~~ | ~~**1093**~~ | ~~vp~~ | | | | |
| 056 | **337** | 871 | | pv | | | ~~020~~ | ~~121~~ | ~~**1087**~~ | ~~vp~~ | | | | |
| 057 | 343 | 865 | | | | vv | ~~024~~ | ~~145~~ | ~~**1063**~~ | ~~vp~~ | | | | |
| 058 | **349** | **859** | pp | | | | ~~028~~ | ~~169~~ | ~~**1039**~~ | ~~vp~~ | | | | |
| 059 | 355 | **853** | | | vp | | ~~029~~ | ~~175~~ | ~~**1033**~~ | ~~vp~~ | | | | |
| 060 | 361 | 847 | | | | vv | ~~031~~ | ~~187~~ | ~~**1021**~~ | ~~vp~~ | | | | |
| 061 | **367** | 841 | | pv | | | ~~036~~ | ~~217~~ | ~~**991**~~ | ~~vp~~ | | | | |
| 062 | **373** | 835 | | pv | | | ~~048~~ | ~~289~~ | ~~**919**~~ | ~~vp~~ | | | | |
| 063 | **379** | **829** | pp | | | | ~~050~~ | ~~301~~ | ~~**907**~~ | ~~vp~~ | | | | |
| 064 | 385 | **823** | | | vp | | ~~054~~ | ~~325~~ | ~~**883**~~ | ~~vp~~ | | | | |
| 065 | 391 | 817 | | | | vv | ~~059~~ | ~~355~~ | ~~**853**~~ | ~~vp~~ | | | | |
| 066 | **397** | **811** | pp | | | | ~~064~~ | ~~385~~ | ~~**823**~~ | ~~vp~~ | | | | |
| 067 | 403 | 805 | | | | vv | ~~075~~ | ~~451~~ | ~~**757**~~ | ~~vp~~ | | | | |
| 068 | **409** | 799 | | pv | | | ~~078~~ | ~~469~~ | ~~**739**~~ | ~~vp~~ | | | | |
| 069 | 415 | 793 | | | | vv | ~~079~~ | ~~475~~ | ~~**733**~~ | ~~vp~~ | | | | |
| 070 | **421** | **787** | pp | | | | ~~080~~ | ~~481~~ | ~~**727**~~ | ~~vp~~ | | | | |
| 071 | 427 | 781 | | | | vv | ~~086~~ | ~~517~~ | ~~**691**~~ | ~~vp~~ | | | | |
| 072 | **433** | 775 | | pv | | | ~~089~~ | ~~535~~ | ~~**673**~~ | ~~vp~~ | | | | |
| 073 | **439** | **769** | pp | | | | ~~094~~ | ~~565~~ | ~~**643**~~ | ~~vp~~ | | | | |
| 074 | 445 | 763 | | | | vv | ~~098~~ | ~~589~~ | ~~**619**~~ | ~~vp~~ | | | | |
| 075 | 451 | **757** | | | vp | | ~~099~~ | ~~595~~ | ~~**613**~~ | ~~vp~~ | | | | |
| 076 | **457** | **751** | pp | | | | 004 | 25 | 1183 | vv | 004 | 25 | 1183 | vv |
| 077 | **463** | 745 | | pv | | | 008 | 49 | 1159 | vv | 008 | 49 | 1159 | vv |
| 078 | 469 | **739** | | | vp | | 022 | 133 | 1075 | vv | 022 | 133 | 1075 | vv |
| 079 | 475 | **733** | | | vp | | 034 | 205 | 1003 | vv | 034 | 205 | 1003 | vv |
| 080 | 481 | **727** | | | vp | | 039 | 235 | 973 | vv | 039 | 235 | 973 | vv |
| 081 | **487** | 721 | | pv | | | 041 | 247 | 961 | vv | 041 | 247 | 961 | vv |
| 082 | 493 | 715 | | | | vv | 042 | 253 | 955 | vv | 042 | 253 | 955 | vv |
| 083 | **499** | **709** | pp | | | | 043 | 259 | 949 | vv | 043 | 259 | 949 | vv |
| 084 | 505 | 703 | | | | vv | 044 | 265 | 943 | vv | 044 | 265 | 943 | vv |
| 085 | 511 | 697 | | | | vv | 049 | 295 | 913 | vv | 049 | 295 | 913 | vv |
| 086 | 517 | **691** | | | vp | | 053 | 319 | 889 | vv | 053 | 319 | 889 | vv |
| 087 | **523** | 685 | | pv | | | 057 | 343 | 865 | vv | 057 | 343 | 865 | vv |
| 088 | 529 | 679 | | | | vv | 060 | 361 | 847 | vv | 060 | 361 | 847 | vv |
| 089 | 535 | **673** | | | vp | | 065 | 391 | 817 | vv | 065 | 391 | 817 | vv |
| 090 | **541** | 667 | | pv | | | 067 | 403 | 805 | vv | 067 | 403 | 805 | vv |
| 091 | **547** | **661** | pp | | | | 069 | 415 | 793 | vv | 069 | 415 | 793 | vv |
| 092 | 553 | 655 | | | | vv | 071 | 427 | 781 | vv | 071 | 427 | 781 | vv |
| 093 | 559 | 649 | | | | vv | 074 | 445 | 763 | vv | 074 | 445 | 763 | vv |
| 094 | 565 | **643** | | | vp | | 082 | 493 | 715 | vv | 082 | 493 | 715 | vv |
| 095 | **571** | 637 | | pv | | | 084 | 505 | 703 | vv | 084 | 505 | 703 | vv |

| | | | | | | | | | | | | | |
|---|---|---|---|---|---|---|---|---|---|---|---|---|---|
| 096 | **577** | **631** | pp | | | 085 | 511 | 697 | ~~vv~~ | ~~085~~ | ~~511~~ | ~~697~~ | ~~vv~~ |
| 097 | 583 | 625 | | | vv | 088 | 529 | 679 | ~~vv~~ | ~~088~~ | ~~529~~ | ~~679~~ | ~~vv~~ |
| 098 | 589 | **619** | | vp | | 092 | 553 | 655 | ~~vv~~ | ~~092~~ | ~~553~~ | ~~655~~ | ~~vv~~ |
| 099 | 595 | **613** | | vp | | 093 | 559 | 649 | ~~vv~~ | ~~093~~ | ~~559~~ | ~~639~~ | ~~vv~~ |
| 100 | **601** | **607** | pp | | | 097 | 583 | 625 | ~~vv~~ | ~~097~~ | ~~583~~ | ~~625~~ | ~~vv~~ |

½#pp = 20
½#(p+v) = 55
½#vv = 25

A = 200, #p = 95, #v = 105, #pp = 40, #(pv+vp) = 110, #vv = 50.

Lösungen zur Zahl 1208
($H_+(200)^{+1}$, $H_+(200)^{-1}$), Zählung, Beispiel zu #v > #p. (Die Fälle 101 bis 200 sind identisch mit denen von 1 bis 100.)
Tabelle 11.2.

Die Darstellung wie in Tabelle 11.2 kann in jedem Fall '#v > #p' durchgeführt werden.

Betrachtung der unsymmetrischen Quadranten 2 und 4.
Zunächst 2. Quadrant. Die Vorgehensweise ist wie oben. Betrachtet wird eine beliebige Rhombusseite ($H_-(A)^{+1}$, $H_+(A)^{-1}$) aus der Lösungsmenge von (2GH). Es gilt:
(22) $0 < \#p_- < A$, denn beide Abschnitte beginnen mit einer Primzahl, bzw. beginnen mit dem Übergewicht der Primzahlen - und die Abschnitte bestehen nicht nur aus Primzahlen.
(23) $0 < \#p_+ < A$, Begründung wie bei (22).

(24) $0 < \#v_- < A$, Begründung wie bei (22).
(25) $0 < \#v_+ < A$, Begründung wie bei (22).

(26) $A = \#p_- + \#v_-$, denn jedes Glied ist entweder Primzahl oder Vielfaches.
(27) $A = \#p_+ + \#v_+$, Begründung wie bei (26).

Bei den Paaren der Rhombusseite ($H_-(A)^{+1}$, $H_+(A)^{-1}$) gibt es nur die Kombinationen $p_-p_+$, $p_-v_+$, $v_-p_+$ und $v_-v_+$ und damit die Zählfunktionen $\#p_-p_+$, $\#p_-v_+$, $\#v_-p_+$ und $\#v_-v_+$. Es ist detaillierter als oben vorauszusetzen:
(28) $\#p_-v_+ \neq \#v_-p_+$
(29) $\#p_- < \#v_-$, $\#p_+ < \#v_+$.

Aus (29) folgt nach dem Schubfachprinzip, daß es $v_-v_+$-Fälle gibt, so daß $\#v_-v_+$ nicht verschwindet. Ausserdem ist $\#v_-v_+$ kleiner als A, weil sonst $\#p_-p_+$ verschwinden würde, was durch (22) ausgeschlossen ist. Es gilt sogar wegen des Schubfachprinzips:
(30) $0 < \#v_- - \#p_- \leq \#v_-v_+ < A$.
(31) $0 < \#v_+ - \#p_+ \leq \#v_-v_+ < A$.
Da die Ausdrücke ($\#v_- - \#p_-$) und ($\#v_+ - \#p_+$) beide kleiner oder gleich $\#v_-v_+$ sind und da
$\#v_- + \#p_- = \#v_+ + \#p_+$, d.h. $\#v_- - \#p_+ = \#v_+ - \#p_-$
gilt auch:
(32) $0 < \#v_- - \#p_+ \leq \#v_-v_+ < A$.
(33) $0 < \#v_+ - \#p_- \leq \#v_-v_+ < A$.

Ein Verschwinden von $\#p_-v_+$ würde bedeuten, daß alle $p_-$ sich mit $p_+$ verbinden - und dann auch alle $v_-$ mit $v_+$. Das kann nicht sein wegen der Unterschiedlichkeit der Komponenten-Abschnitte. Ebenso wenig kann $\#p_-v_+$ gleich A sein. Das gleiche gilt für $\#v_-p_+$:
(34) $0 < \#p_-v_+ < A$.
(35) $0 < \#v_-p_+ < A$.

Ebenso kann nicht $\#p_-p_+$ gleich A sein. Denn dann gäbe es in dem Komponenten-Abschnitt keine Vielfachen. Ob allerdings $\#p_-p_+$ verschwinden kann ist offen - und ist die Frage von GOLDBACH:
(36) $0 \leq \#p_-p_+ < A$.

Die Zählfunktionen liegen also alle zwischen 0 und A - nur bei $\#p_-p_+$ ist das Verschwinden noch nicht ausgeschlossen und ist Gegenstand der Betrachtung.

Allgemein gilt:
(37) $\#v_-v_+ + \#v_-p_+ = \#v_-$ und $\#p_-p_+ + \#p_-v_+ = \#p_-$, d.h. $\#v_- - \#v_-v_+ = \#v_-p_+$, $\#p_- - \#p_-p_+ = \#p_-v_+$
(38) $\#v_-v_+ + \#p_-v_+ = \#v_+$ und $\#p_-p_+ + \#v_-p_+ = \#p_+$, d.h. $\#v_+ - \#v_-v_+ = \#p_-v_+$, $\#p_+ - \#p_-p_+ = \#v_-p_+$

Daraus folgt:
(39) $\#v_- - \#p_+ = \#v_-v_+ - \#p_-p_+$
(40) $\#v_+ - \#p_- = \#v_-v_+ - \#p_-p_+$

Die Betrachtung von (35), (33), (39) und (40) zeigt: nur wenn $\#v_- - \#p_+$ oder $\#v_+ - \#p_- = \#v_-v_+$ ist, kann $\#p_-p_+ = 0$ sein. Ist aber $\#v_- - \#p_+$ oder $\#v_+ - \#p_- < \#v_-v_+$, dann ist zwangsläufig $\#p_-p_+ > 0$. Zur Bestätigung der (2GH) ist also
(41) $\#v_- - \#p_+ < \#v_-v_+$,
(42) $\#v_+ - \#p_- < \#v_-v_+$
zu beweisen.

Zum allgemeinen Beweis wird $\#v_-v_+$ mit Hilfe der relativen Häufigkeiten abgeschätzt. Diese Prozedur ist ebenso wie in (11) bis (23). Dabei werden $H_-$ und $H_+$ hinsichtlich der Primzahl-Verteilung als hinreichend gleich angenommen. Damit ist dann die (2GH) auch im 2. und 4. Quadranten bestätigt.

Ferner gilt:
(43) $\#p_-p_+ = [A - (\#p_-v_+ + \#v_- + p_+) - (\#v_- - \#p_+)] / 2$.
Und außerdem ist dann:
(44) $\#v_-v_+ = [A - (\#p_-v_+ + \#v_- + p_+) - (\#v_- - \#p_+)] / 2 + [\#v_- - \#p_+]$.
Beispiel A = 191:

Nr $(H_+(191)^{+1}, H_+(191)^{-1})$

Originale

| Nr | | | pp | pv | vp | vv |
|---|---|---|---|---|---|---|
| 001 | **5** | 1147 | | pv | | |
| 002 | **11** | 1141 | | pv | | |
| 003 | **17** | 1135 | | pv | | |
| 004 | **23** | **1129** | pp | | | |
| 005 | **29** | **1123** | pp | | | |
| 006 | 35 | **1117** | | | vp | |
| 007 | **41** | 1111 | | pv | | |
| 008 | **47** | 1105 | | pv | | |
| 009 | **53** | 1099 | | pv | | |
| 010 | **59** | **1093** | pp | | | |
| 011 | 65 | **1087** | | | vp | |
| 012 | **71** | 1081 | | pv | | |
| 013 | 77 | 1075 | | | | vv |
| 014 | **83** | **1069** | pp | | | |
| 015 | **89** | **1063** | pp | | | |
| 016 | 95 | 1057 | | | | vv |
| 017 | **101** | **1051** | pp | | | |
| 018 | **107** | 1045 | | pv | | |
| 019 | **113** | **1039** | pp | | | |
| 020 | 119 | **1033** | | | vp | |
| 021 | 125 | 1027 | | | | vv |
| 022 | **131** | **1021** | pp | | | |
| 023 | **137** | 1015 | | pv | | |
| 024 | 143 | **1009** | | | vp | |
| 025 | **149** | 1003 | | pv | | |
| 026 | 155 | **997** | | | vp | |
| 027 | 161 | **991** | | | vp | |

Nr $(H_+(191)^{+1}, H_+(191)^{-1})$

sortiert nach Typ
pv und vp gestrichen

| Nr | | | pp | pv | vp | vv |
|---|---|---|---|---|---|---|
| 004 | **23** | **1129** | pp | | | |
| 005 | **29** | **1123** | pp | | | |
| 010 | **59** | **1093** | pp | | | |
| 014 | **83** | **1069** | pp | | | |
| 015 | **89** | **1063** | pp | | | |
| 017 | **101** | **1051** | pp | | | |
| 019 | **113** | **1039** | pp | | | |
| 022 | **131** | **1021** | pp | | | |
| 039 | **233** | **919** | pp | | | |
| 045 | **269** | **883** | pp | | | |
| 049 | **293** | **859** | pp | | | |
| 064 | **383** | **769** | pp | | | |
| 067 | **401** | **751** | pp | | | |
| 070 | **419** | **733** | pp | | | |
| 074 | **443** | **709** | pp | | | |
| 077 | **461** | **691** | pp | | | |
| 080 | **479** | **673** | pp | | | |
| 082 | **491** | **661** | pp | | | |
| 085 | **509** | **643** | pp | | | |
| 087 | **521** | **631** | pp | | | |
| 109 | **653** | **499** | pp | | | |
| 120 | **719** | **433** | pp | | | |
| 124 | **743** | **409** | pp | | | |
| 129 | **773** | **379** | pp | | | |
| 137 | **821** | **331** | pp | | | |
| 140 | **839** | **313** | pp | | | |
| 147 | **881** | **271** | pp | | | |

Nr $(H_+(191)^{+1}, H_+(191)^{-1})$

sortiert nach Typ
ohne pv und vp und
Überstand vv gestrichen

| Nr | | | pp | pv | vp | vv |
|---|---|---|---|---|---|---|
| 004 | **23** | **1129** | pp | | | |
| 005 | **29** | **1123** | pp | | | |
| 010 | **59** | **1093** | pp | | | |
| 014 | **83** | **1069** | pp | | | |
| 015 | **89** | **1063** | pp | | | |
| 017 | **101** | **1051** | pp | | | |
| 019 | **113** | **1039** | pp | | | |
| 022 | **131** | **1021** | pp | | | |
| 039 | **233** | **919** | pp | | | |
| 045 | **269** | **883** | pp | | | |
| 049 | **293** | **859** | pp | | | |
| 064 | **383** | **769** | pp | | | |
| 067 | **401** | **751** | pp | | | |
| 070 | **419** | **733** | pp | | | |
| 074 | **443** | **709** | pp | | | |
| 077 | **461** | **691** | pp | | | |
| 080 | **479** | **673** | pp | | | |
| 082 | **491** | **661** | pp | | | |
| 085 | **509** | **643** | pp | | | |
| 087 | **521** | **631** | pp | | | |
| 109 | **653** | **499** | pp | | | |
| 120 | **719** | **433** | pp | | | |
| 124 | **743** | **409** | pp | | | |
| 129 | **773** | **379** | pp | | | |
| 137 | **821** | **331** | pp | | | |
| 140 | **839** | **313** | pp | | | |
| 147 | **881** | **271** | pp | | | |

| | | | |
|---|---|---|---|
| 028 | **167** | 985 | pv |
| 029 | **173** | 979 | pv |
| 030 | **179** | 973 | pv |
| 031 | 185 | **967** | vp |
| 032 | **191** | 961 | pv |
| 033 | **197** | 955 | pv |
| 034 | 203 | 949 | vv |
| 035 | 209 | 943 | vv |
| 036 | 215 | **937** | vp |
| 037 | 221 | 931 | vv |
| 038 | **227** | 925 | pv |
| 039 | **233** | **919** | pp |
| 040 | **239** | 913 | pv |
| 041 | 245 | **907** | vp |
| 042 | **251** | 901 | pv |
| 043 | **257** | 895 | pv |
| 044 | **263** | 889 | pv |
| 045 | **269** | **883** | pp |
| 046 | 275 | **877** | vp |
| 047 | **281** | 871 | pv |
| 048 | 287 | 865 | vv |
| 049 | **293** | **859** | pp |
| 050 | 299 | **853** | vp |
| 051 | 305 | 847 | vv |
| 052 | **311** | 841 | pv |
| 053 | **317** | 835 | pv |
| 054 | 323 | **829** | vp |
| 055 | 329 | **823** | vp |
| 056 | 335 | 817 | vv |
| 057 | 341 | **811** | vp |
| 058 | **347** | 805 | pv |
| 059 | **353** | 799 | pv |
| 060 | **359** | 793 | pv |
| 061 | 365 | **787** | vp |
| 062 | 371 | 781 | vv |
| 063 | 377 | 775 | vv |
| 064 | **383** | **769** | pp |
| 065 | **389** | 763 | pv |
| 066 | 395 | **757** | vp |
| 067 | **401** | **751** | pp |
| 068 | 407 | 745 | vv |
| 069 | 413 | **739** | vp |
| 070 | **419** | **733** | pp |
| 071 | 425 | **727** | vp |
| 072 | **431** | 721 | pv |
| 073 | 437 | 715 | vv |
| 074 | **443** | **709** | pp |
| 075 | **449** | 703 | pv |
| 076 | 455 | 697 | vv |
| 077 | **461** | **691** | pp |
| 078 | **467** | 685 | pv |
| 079 | 473 | 679 | vv |
| 080 | **479** | **673** | pp |
| 081 | 485 | 667 | vv |
| 082 | **491** | **661** | pp |
| 083 | 497 | 655 | vv |
| 084 | **503** | 649 | pv |
| 085 | **509** | **643** | pp |
| 086 | 515 | 637 | vv |
| 087 | **521** | **631** | pp |
| 088 | 527 | 625 | vv |
| 089 | 533 | **619** | vp |
| 090 | 539 | **613** | vp |
| 091 | 545 | **607** | vp |
| 092 | 551 | **601** | vp |
| 093 | **557** | 595 | pv |
| 094 | **563** | 589 | pv |
| 095 | **569** | 583 | pv |
| 096 | 575 | **577** | vp |
| 097 | 581 | **571** | vp |
| 098 | **587** | 565 | pv |
| 152 | **911** | **241** | pp |
| 155 | **929** | **223** | pp |
| 157 | **941** | **211** | pp |
| 159 | **953** | **199** | pp |
| 162 | **971** | **181** | pp |
| 169 | **1013** | **139** | pp |
| 175 | **1049** | **103** | pp |
| 182 | **1091** | **61** | pp |
| 185 | **1109** | **43** | pp |
| ~~001~~ | ~~**5**~~ | ~~1147~~ | ~~pv~~ |
| ~~002~~ | ~~**11**~~ | ~~1141~~ | ~~pv~~ |
| ~~003~~ | ~~**17**~~ | ~~1135~~ | ~~pv~~ |
| ~~007~~ | ~~**41**~~ | ~~1111~~ | ~~pv~~ |
| ~~008~~ | ~~**47**~~ | ~~1105~~ | ~~pv~~ |
| ~~009~~ | ~~**53**~~ | ~~1099~~ | ~~pv~~ |
| ~~012~~ | ~~**71**~~ | ~~1081~~ | ~~pv~~ |
| ~~018~~ | ~~**107**~~ | ~~1045~~ | ~~pv~~ |
| ~~023~~ | ~~**137**~~ | ~~1015~~ | ~~pv~~ |
| ~~025~~ | ~~**149**~~ | ~~1003~~ | ~~pv~~ |
| ~~028~~ | ~~**167**~~ | ~~985~~ | ~~pv~~ |
| ~~029~~ | ~~**173**~~ | ~~979~~ | ~~pv~~ |
| ~~030~~ | ~~**179**~~ | ~~973~~ | ~~pv~~ |
| ~~032~~ | ~~**191**~~ | ~~961~~ | ~~pv~~ |
| ~~033~~ | ~~**197**~~ | ~~955~~ | ~~pv~~ |
| ~~038~~ | ~~**227**~~ | ~~925~~ | ~~pv~~ |
| ~~041~~ | ~~**239**~~ | ~~913~~ | ~~pv~~ |
| ~~042~~ | ~~**251**~~ | ~~901~~ | ~~pv~~ |
| ~~043~~ | ~~**257**~~ | ~~895~~ | ~~pv~~ |
| ~~044~~ | ~~**263**~~ | ~~889~~ | ~~pv~~ |
| ~~047~~ | ~~**281**~~ | ~~871~~ | ~~pv~~ |
| ~~052~~ | ~~**311**~~ | ~~841~~ | ~~pv~~ |
| ~~053~~ | ~~**317**~~ | ~~835~~ | ~~pv~~ |
| ~~058~~ | ~~**347**~~ | ~~805~~ | ~~pv~~ |
| ~~059~~ | ~~**353**~~ | ~~799~~ | ~~pv~~ |
| ~~060~~ | ~~**359**~~ | ~~793~~ | ~~pv~~ |
| ~~065~~ | ~~**389**~~ | ~~763~~ | ~~pv~~ |
| ~~072~~ | ~~**431**~~ | ~~721~~ | ~~pv~~ |
| ~~075~~ | ~~**449**~~ | ~~703~~ | ~~pv~~ |
| ~~078~~ | ~~**467**~~ | ~~685~~ | ~~pv~~ |
| ~~084~~ | ~~**503**~~ | ~~649~~ | ~~pv~~ |
| ~~093~~ | ~~**557**~~ | ~~595~~ | ~~pv~~ |
| ~~094~~ | ~~**563**~~ | ~~589~~ | ~~pv~~ |
| ~~095~~ | ~~**569**~~ | ~~583~~ | ~~pv~~ |
| ~~098~~ | ~~**587**~~ | ~~565~~ | ~~pv~~ |
| ~~099~~ | ~~**593**~~ | ~~559~~ | ~~pv~~ |
| ~~100~~ | ~~**599**~~ | ~~553~~ | ~~pv~~ |
| ~~103~~ | ~~**617**~~ | ~~535~~ | ~~pv~~ |
| ~~107~~ | ~~**641**~~ | ~~511~~ | ~~pv~~ |
| ~~108~~ | ~~**647**~~ | ~~505~~ | ~~pv~~ |
| ~~110~~ | ~~**659**~~ | ~~493~~ | ~~pv~~ |
| ~~113~~ | ~~**677**~~ | ~~475~~ | ~~pv~~ |
| ~~114~~ | ~~**683**~~ | ~~469~~ | ~~pv~~ |
| ~~117~~ | ~~**701**~~ | ~~451~~ | ~~pv~~ |
| ~~127~~ | ~~**761**~~ | ~~391~~ | ~~pv~~ |
| ~~133~~ | ~~**797**~~ | ~~355~~ | ~~pv~~ |
| ~~135~~ | ~~**809**~~ | ~~343~~ | ~~pv~~ |
| ~~138~~ | ~~**827**~~ | ~~325~~ | ~~pv~~ |
| ~~143~~ | ~~**857**~~ | ~~295~~ | ~~pv~~ |
| ~~144~~ | ~~**863**~~ | ~~289~~ | ~~pv~~ |
| ~~148~~ | ~~**887**~~ | ~~265~~ | ~~pv~~ |
| ~~158~~ | ~~**947**~~ | ~~205~~ | ~~pv~~ |
| ~~163~~ | ~~**977**~~ | ~~175~~ | ~~pv~~ |
| ~~164~~ | ~~**983**~~ | ~~169~~ | ~~pv~~ |
| ~~170~~ | ~~**1019**~~ | ~~133~~ | ~~pv~~ |
| ~~172~~ | ~~**1031**~~ | ~~121~~ | ~~pv~~ |
| ~~177~~ | ~~**1061**~~ | ~~91~~ | ~~pv~~ |
| ~~183~~ | ~~**1097**~~ | ~~55~~ | ~~pv~~ |
| ~~184~~ | ~~**1103**~~ | ~~49~~ | ~~pv~~ |
| ~~006~~ | ~~35~~ | ~~**1117**~~ | ~~vp~~ |
| ~~011~~ | ~~65~~ | ~~**1087**~~ | ~~vp~~ |
| ~~020~~ | ~~119~~ | ~~**1033**~~ | ~~vp~~ |
| 152 | **911** | **241** | pp |
| 155 | **929** | **223** | pp |
| 157 | **941** | **211** | pp |
| 159 | **953** | **199** | pp |
| 162 | **971** | **181** | pp |
| 169 | **1013** | **139** | pp |
| 175 | **1049** | **103** | pp |
| 182 | **1091** | **61** | pp |
| 185 | **1109** | **43** | pp |

| | | | | | | | | | | | | | | |
|---|---|---|---|---|---|---|---|---|---|---|---|---|---|---|
| 099 | **593** | 559 | | pv | | | ~~024~~ | ~~143~~ | ~~**1009**~~ | ~~vp~~ | | | | |
| 100 | **599** | 553 | | pv | | | ~~026~~ | ~~155~~ | ~~**997**~~ | ~~vp~~ | | | | |
| 101 | 605 | **547** | | | vp | | ~~027~~ | ~~161~~ | ~~**991**~~ | ~~vp~~ | | | | |
| 102 | 611 | **541** | | | vp | | ~~031~~ | ~~185~~ | ~~**967**~~ | ~~vp~~ | | | | |
| 103 | **617** | 535 | | pv | | | ~~036~~ | ~~215~~ | ~~**937**~~ | ~~vp~~ | | | | |
| 104 | 623 | 529 | | | | vv | ~~041~~ | ~~245~~ | ~~**907**~~ | ~~vp~~ | | | | |
| 105 | 629 | **523** | | | vp | | ~~046~~ | ~~275~~ | ~~**877**~~ | ~~vp~~ | | | | |
| 106 | 635 | 517 | | | | vv | ~~050~~ | ~~299~~ | ~~**853**~~ | ~~vp~~ | | | | |
| 107 | **641** | 511 | | pv | | | ~~054~~ | ~~323~~ | ~~**829**~~ | ~~vp~~ | | | | |
| 108 | **647** | 505 | | pv | | | ~~055~~ | ~~329~~ | ~~**823**~~ | ~~vp~~ | | | | |
| 109 | **653** | **499** | pp | | | | ~~057~~ | ~~341~~ | ~~**811**~~ | ~~vp~~ | | | | |
| 110 | **659** | 493 | | pv | | | ~~061~~ | ~~365~~ | ~~**787**~~ | ~~vp~~ | | | | |
| 111 | 665 | **487** | | | vp | | ~~066~~ | ~~395~~ | ~~**757**~~ | ~~vp~~ | | | | |
| 112 | 671 | 481 | | | | vv | ~~069~~ | ~~413~~ | ~~**739**~~ | ~~vp~~ | | | | |
| 113 | **677** | 475 | | pv | | | ~~071~~ | ~~425~~ | ~~**727**~~ | ~~vp~~ | | | | |
| 114 | **683** | 469 | | pv | | | ~~089~~ | ~~533~~ | ~~**619**~~ | ~~vp~~ | | | | |
| 115 | 689 | **463** | | | vp | | ~~090~~ | ~~539~~ | ~~**613**~~ | ~~vp~~ | | | | |
| 116 | 695 | **457** | | | vp | | ~~091~~ | ~~545~~ | ~~**607**~~ | ~~vp~~ | | | | |
| 117 | **701** | 451 | | pv | | | ~~092~~ | ~~551~~ | ~~**601**~~ | ~~vp~~ | | | | |
| 118 | 707 | 445 | | | | vv | ~~096~~ | ~~575~~ | ~~**577**~~ | ~~vp~~ | | | | |
| 119 | 713 | **439** | | | vp | | ~~097~~ | ~~581~~ | ~~**571**~~ | ~~vp~~ | | | | |
| 120 | **719** | **433** | pp | | | | ~~101~~ | ~~605~~ | ~~**547**~~ | ~~vp~~ | | | | |
| 121 | 725 | 427 | | | | vv | ~~102~~ | ~~611~~ | ~~**541**~~ | ~~vp~~ | | | | |
| 122 | 731 | **421** | | | vp | | ~~105~~ | ~~629~~ | ~~**523**~~ | ~~vp~~ | | | | |
| 123 | 737 | 415 | | | | vv | ~~111~~ | ~~665~~ | ~~**487**~~ | ~~vp~~ | | | | |
| 124 | **743** | **409** | pp | | | | ~~115~~ | ~~689~~ | ~~**463**~~ | ~~vp~~ | | | | |
| 125 | 749 | 403 | | | | vv | ~~116~~ | ~~695~~ | ~~**457**~~ | ~~vp~~ | | | | |
| 126 | 755 | **397** | | | vp | | ~~119~~ | ~~713~~ | ~~**439**~~ | ~~vp~~ | | | | |
| 127 | **761** | 391 | | pv | | | ~~122~~ | ~~731~~ | ~~**421**~~ | ~~vp~~ | | | | |
| 128 | 767 | 385 | | | | vv | ~~126~~ | ~~755~~ | ~~**397**~~ | ~~vp~~ | | | | |
| 129 | **773** | **379** | pp | | | | ~~130~~ | ~~779~~ | ~~**373**~~ | ~~vp~~ | | | | |
| 130 | 779 | **373** | | | vp | | ~~131~~ | ~~785~~ | ~~**367**~~ | ~~vp~~ | | | | |
| 131 | 785 | **367** | | | vp | | ~~134~~ | ~~803~~ | ~~**349**~~ | ~~vp~~ | | | | |
| 132 | 791 | 361 | | | | vv | ~~136~~ | ~~815~~ | ~~**337**~~ | ~~vp~~ | | | | |
| 133 | **797** | 355 | | pv | | | ~~141~~ | ~~845~~ | ~~**307**~~ | ~~vp~~ | | | | |
| 134 | 803 | **349** | | | vp | | ~~145~~ | ~~869~~ | ~~**283**~~ | ~~vp~~ | | | | |
| 135 | **809** | 343 | | pv | | | ~~146~~ | ~~875~~ | ~~**277**~~ | ~~vp~~ | | | | |
| 136 | 815 | **337** | | | vp | | ~~154~~ | ~~823~~ | ~~**229**~~ | ~~vp~~ | | | | |
| 137 | **821** | **331** | pp | | | | ~~160~~ | ~~859~~ | ~~**193**~~ | ~~vp~~ | | | | |
| 138 | **827** | 325 | | pv | | | ~~165~~ | ~~989~~ | ~~**163**~~ | ~~vp~~ | | | | |
| 139 | 833 | 319 | | | | vv | ~~166~~ | ~~995~~ | ~~**157**~~ | ~~vp~~ | | | | |
| 140 | **839** | **313** | pp | | | | ~~167~~ | ~~1001~~ | ~~**151**~~ | ~~vp~~ | | | | |
| 141 | 845 | **307** | | | vp | | ~~171~~ | ~~1025~~ | ~~**127**~~ | ~~vp~~ | | | | |
| 142 | 851 | 301 | | | | vv | ~~174~~ | ~~1043~~ | ~~**109**~~ | ~~vp~~ | | | | |
| 143 | **857** | 295 | | pv | | | ~~176~~ | ~~1055~~ | ~~**97**~~ | ~~vp~~ | | | | |
| 144 | **863** | 289 | | pv | | | ~~179~~ | ~~1073~~ | ~~**79**~~ | ~~vp~~ | | | | |
| 145 | 869 | **283** | | | vp | | ~~180~~ | ~~1079~~ | ~~**73**~~ | ~~vp~~ | | | | |
| 146 | 875 | **277** | | | vp | | ~~181~~ | ~~1085~~ | ~~**67**~~ | ~~vp~~ | | | | |
| 147 | **881** | **271** | pp | | | | ~~186~~ | ~~1115~~ | ~~**37**~~ | ~~vp~~ | | | | |
| 148 | **887** | 265 | | pv | | | ~~187~~ | ~~1121~~ | ~~**31**~~ | ~~vp~~ | | | | |
| 149 | 893 | 259 | | | | vv | ~~189~~ | ~~1133~~ | ~~**19**~~ | ~~vp~~ | | | | |
| 150 | 899 | 253 | | | | vv | ~~190~~ | ~~1139~~ | ~~**13**~~ | ~~vp~~ | | | | |
| 151 | 905 | 247 | | | | vv | ~~191~~ | ~~1145~~ | ~~**7**~~ | ~~vp~~ | | | | |
| 152 | **911** | **241** | pp | | | | 013 | 77 | 1075 | vv | 013 | 77 | 1075 | vv |
| 153 | 917 | 235 | | | | vv | 016 | 95 | 1057 | vv | 016 | 95 | 1057 | vv |
| 154 | 923 | **229** | | | vp | | 021 | 125 | 1027 | vv | 021 | 125 | 1027 | vv |
| 155 | **929** | **223** | pp | | | | 034 | 203 | 949 | vv | 034 | 203 | 949 | vv |
| 156 | 935 | 217 | | | | vv | 035 | 209 | 943 | vv | 035 | 209 | 943 | vv |
| 157 | **941** | **211** | pp | | | | 037 | 221 | 931 | vv | 037 | 221 | 931 | vv |
| 158 | **947** | 205 | | pv | | | 048 | 287 | 865 | vv | 048 | 287 | 865 | vv |
| 159 | **953** | **199** | pp | | | | 051 | 305 | 847 | vv | 051 | 305 | 847 | vv |
| 160 | 959 | **193** | | | vp | | 056 | 335 | 817 | vv | 056 | 335 | 817 | vv |
| 161 | 965 | 187 | | | | vv | 062 | 371 | 781 | vv | 062 | 371 | 781 | vv |
| 162 | **971** | **181** | pp | | | | 063 | 377 | 775 | vv | 063 | 377 | 775 | vv |
| 163 | **977** | 175 | | pv | | | 068 | 407 | 745 | vv | 068 | 407 | 745 | vv |
| 164 | **983** | 169 | | pv | | | 073 | 437 | 715 | vv | 073 | 437 | 715 | vv |
| 165 | 989 | **163** | | | vp | | 076 | 455 | 697 | vv | 076 | 455 | 697 | vv |
| 166 | 995 | **157** | | | vp | | 079 | 473 | 679 | vv | 079 | 473 | 679 | vv |
| 167 | 1001 | **151** | | | vp | | 081 | 485 | 667 | vv | 081 | 485 | 667 | vv |
| 168 | 1007 | 145 | | | | vv | 083 | 497 | 655 | vv | 083 | 497 | 655 | vv |
| 169 | **1013** | **139** | pp | | | | 086 | 515 | 637 | vv | 086 | 515 | 637 | vv |

| | | | | | | | | | | | | | | | |
|---|---|---|---|---|---|---|---|---|---|---|---|---|---|---|---|
| 170 | **1019** | 133 | | pv | | | 088 | 527 | 625 | vv | 088 | 527 | 625 | vv |
| 171 | 1025 | **127** | | | vp | | 104 | 623 | 529 | vv | 104 | 623 | 529 | vv |
| 172 | **1031** | 121 | | pv | | | 106 | 635 | 517 | vv | 106 | 635 | 517 | vv |
| 173 | 1037 | 115 | | | | vv | 112 | 671 | 481 | vv | 112 | 671 | 481 | vv |
| 174 | 1043 | **109** | | | vp | | 118 | 707 | 445 | vv | 118 | 707 | 445 | vv |
| 175 | **1049** | **103** | pp | | | | 121 | 725 | 427 | vv | 121 | 725 | 427 | vv |
| 176 | 1055 | **97** | | | vp | | 123 | 737 | 415 | vv | 123 | 737 | 415 | vv |
| 177 | **1061** | 91 | pv | | | | 125 | 749 | 403 | vv | 125 | 749 | 403 | vv |
| 178 | 1067 | 85 | | | | vv | 128 | 767 | 385 | vv | 128 | 767 | 385 | vv |
| 179 | 1073 | **79** | | | vp | | 132 | 791 | 361 | vv | 132 | 791 | 361 | vv |
| 180 | 1079 | **73** | | | vp | | 139 | 833 | 319 | vv | 139 | 833 | 319 | vv |
| 181 | 1085 | **67** | | | vp | | 142 | 851 | 301 | vv | 142 | 851 | 301 | vv |
| 182 | **1091** | **61** | pp | | | | 149 | 893 | 259 | vv | 149 | 893 | 259 | vv |
| 183 | **1097** | 55 | | pv | | | 150 | 899 | 253 | vv | 150 | 899 | 253 | vv |
| 184 | **1103** | 49 | | pv | | | 151 | 905 | 247 | vv | 151 | 905 | 247 | vv |
| 185 | **1109** | **43** | pp | | | | 153 | 917 | 235 | vv | 153 | 917 | 235 | vv |
| 186 | 1115 | **37** | | | vp | | 156 | 935 | 217 | vv | 56 | 935 | 217 | vv |
| 187 | 1121 | **31** | | | vp | | 161 | 965 | 187 | vv | 161 | 965 | 187 | vv |
| 188 | 1127 | 25 | | | | vv | ~~168~~ | ~~1007~~ | ~~145~~ | ~~vv~~ | | | | |
| 189 | 1133 | **19** | | | vp | | ~~173~~ | ~~1037~~ | ~~115~~ | ~~vv~~ | | | | |
| 190 | 1139 | **13** | | | vp | | ~~178~~ | ~~1067~~ | ~~85~~ | ~~vv~~ | | | | |
| 191 | 1145 | **7** | | | vp | | ~~188~~ | ~~1127~~ | ~~25~~ | ~~vv~~ | | | | |

$A = \#p_+ + \#v_- = \#p_- + \#v_+ = p_+p_- + p_+v_- + v_+p_- + v_+v_-$

$191 = 92 < 99 = 95 < 96 = 36 \quad 59 \quad 56 \quad 40$

Lösungen zur Zahl 1152

($H_-(191)^{+1}$, $H_+(191)^{-1}$), Zählung, Beispiel zu #v > #p. Primzahlen fett.

Tabelle 11.3.

Der Fall #pp = 0 ist als extrem unwahrscheinlicher Fall durch das Zahlenwerk der gegenläufigen Vielfachen-Module der Primzahlen ausgeschlossen.

**Satz 11.2. Bestätigung der Goldbachschen Vermutung.**
*Die Lösungen der (2GH) sind darstellbar als Punkte im Primzahl-Gitter. Der geometrische Ort der Lösungsmenge der (2GH) im Primzahl-Gitter ist ein System von achsen-symmetrischen Rhomben. Es gilt: 1.) Jede gerade Zahl besitzt mindestens eine, im Allgemeinen mehrere (2GH)-Lösungen. 2.) Alle (2GH)-Lösungen, die eine bestimmte gerade Zahl darstellen und nur diese, befinden sich auf einer ihnen umkehrbar eindeutig zugeordneten Rhombus-Seite. 3.) Die Gestalt der (2GH)-Lösungen ist durch ihre Rhombus-Seite bestimmt. 4.) Man kann die Lösungen anhand der Rhombus-Seite explizit angeben. 5.) Es gibt für die zu einer geraden Zahl gehörenden Lösungen - außer ihrer Rhombus-Seite - keine geschlossene Darstellung. – Die (2GH) ist zutreffend.*

## § 12. Ergebnisse: Sätze zu Strukturen in P

### Zu § 1. Die kleinste Obermenge von P

**Definition 1.1.** Primzahl (herkömmliche Definition):
*Primzahl: ≠1, positiv, ganz, nur durch 1 und sich selbst teilbar.*
**Satz 1.1.** Unendlichkeit P:
*Es gibt unendlich viele Primzahlen* (EUKLID).
**Definition 1.2.** Die Folge H:
*H := (±3·2; 1)* = {…, -35, -29, -23, -17,-11, -5, 1, 7, 13, 19, 25, 31, 37, …}.
**Satz 1.2.** H kleinste Obermenge:
*H ist in Z die kleinste aZ, die Obermenge von P \ {3, 2} ist.*
**Satz 1.3.** H Halbgruppe:
*H ist eine kommutative Halbgruppe mit 1-Element.*
**Satz 1.4.** Gestalt in H:
*In H haben alle Glieder die Gestalt* $3 \cdot 2 \cdot k + 1$ *mit* $k \in Z$.
**Definition 1.3.** Primzahl (hier vorgeschlagene Definition):
*Primzahl: Glied von H, nur durch 1 und sich selbst teilbar.*
**Definition 1.4.** Äste von H:
$H_+ := (+3 \cdot 2; 1) = \{1, 7, 13, \dots\}$, $H_- := (+3 \cdot 2; -1) = \{-1, 5, 11, \dots\}$.
**Satz 1.5.** *Gestalt in* $H_+$ *und* $H_-$:
*Für* $x \in H_+$ *ex.* $k \in N$ *mit* $x = 3 \cdot 2 \cdot k + 1$. *Für* $x \in H_-$ *ex.* $k \in N$ *mit* $x = 3 \cdot 2 \cdot k - 1$.
**Satz 1.6.** Primzahlen in der arithmetischen Progression (DIRICHLET, Satz III):
$H_+$ *enthält* ∞ *viele Primzahlen, und zwar die mit der Gestalt* $p = 3 \cdot 2 \cdot k + 1$, $k \in N$.
$H_-$ *enthält* ∞ *viele Primzahlen, und zwar die mit der Gestalt* $p = 3 \cdot 2 \cdot k - 1$, $k \in N$.

**Zu § 2. Das Primzahl-Gitter.**

**Definition 2.1.** *Primzahl-Gitter 2-dimensional* $:= H \times H = \{(x, y) \mid x, y \in H\}$, *…3-dimensional:* $H \times H \times H := \{(x, y, z) \mid x, y, z \in H\}$. *Usw.*

**Satz 2.1.** *Eigenschaften des Primzahl-Gitters* $H^2$*: 1) Sein Koordinaten-Ursprung ist der Punkt (+1, + 1). 2 ) Es wird aufgespannt von H, bzw. von den Ästen von H. 3) Seine Träger-Folgen* $S_-$ *und* $S_+$ *sind einander gleich und sind = N. 4) Es hat die vier Quadranten 1.* $H_+ \times H_+$*, 2.* $H_+ \times H_-$*, 3.* $H_- \times H_-$ *und 4.* $H_- \times H_+$*. 5) Symmetrie besteht nur zwischen dem 1. und 3. Quadranten, d. h. zur Winkelhalbierenden der Quadranten. 6) Es besitzt eine bestimmte geometrische Struktur: Primzahlen und Vielfache generieren Rechtecke und Streifen. - Alles Analog in* $H^3$, $H^n$.

**Zu § 3. Geometrischer Ort der Lösungsmenge der 2-gliedrigen GOLDBACH-Vermutung**.

**Satz 3.1.** *Die Reste -2, 0, +2 (mod 3·2) der geraden Zahlen:* $x \in (2; 10) \Rightarrow s$ *existiert mit* $x := 3 \cdot 2 \cdot s - 2$ *oder* $x := 3 \cdot 2 \cdot s \pm 0$ *oder* $x = :3 \cdot 2 \cdot s + 2$, $1 < s \in N$. *D.h.:* $x \equiv -2$ *oder* $\equiv 0$ *oder* $\equiv +2 \pmod{3 \cdot 2}$.

**Satz 3.2.** *Zerlegung von (+2; +10) in die Folgen* $E_{-2}$, $E_0$ *und* $E_{+2}$*:* $E_{+2} := (3 \cdot 2; 14) = \{14, 20, 26, …\} \equiv +2 \pmod{3 \cdot 2}$, $E_0 := (3 \cdot 2; 12) = \{12, 18, 24, …\} \equiv 0 \pmod{3 \cdot 2}$, $E_{-2} := (3 \cdot 2; 10) = 10, 16, 22, … \equiv -2 \pmod{3 \cdot 2}$. $E_{-2}$, $E_0$ *und* $E_{+2}$ *sind zueinander fremd und vereinigt gleich* $e_n$*. In* $e_n$ *kommt jede gerade Zahl e > 8 genau 1 Mal vor.*

**Satz 3.3.** Die Träger-Folgen $S(E_{+2}) = S(E_0) = S(E_{-2}) = N \setminus \{0, 1\}$.

**Satz 3.4.** Fälle: *Fall 1:* $x, y \in H_+ \Rightarrow x+y \in E_{+2}$. *Fall 2:* $x, y \in H_- \Rightarrow x+y \in E_{-2}$. *Fall 3:* $x \in H_+$, $y \in H_-$*, oder umgekehrt* $\Rightarrow x+y \in E_0$.

**Satz 3.5.** *Die Lösungen von (2GH): Wenn* $p + q = e$, $(3 < p, q) \wedge (p, q \in P) \wedge (e \in E_{+2}$ *oder* $E_{-2}$ *oder* $E_0)$*, so (p, q) Lösung von (2GH).*

**Satz 3.6.** *Gestalt der Lösungen: Fall 1:* $3 \cdot 2 \cdot s_p + 1 + 3 \cdot 2 \cdot s_q + 1 = 3 \cdot 2 \cdot s_e \equiv +2 \pmod{3 \cdot 2}$, *Fall 2:* $3 \cdot 2 \cdot s_p - 1 + 3 \cdot 2 \cdot s_q - 1 = 3 \cdot 2 \cdot s_e \equiv -2$ (mod $\cdot 2$),

**Satz 3.7.** Träger-Schreibweise: *Einheitliche Darstellung der Form der Lösungen mittels der Träger-Schreibweise:* $s_p + s_q = s_e$.

**Definition 3.1.** *Der Komponenten-Abschnitt* $H_{…}(m)^{…}$ *ist ein Anfangs-Abschnitt d.h. Abschnitt der ersten m Glieder von* $H_+$ *oder* $H_-$*. 1) Er ist ein m-tupel. 2) Suffix +: er stammt von* $H_+$*, Suffix -: er stammt von* $H_-$*. 3) Exponent +1: Reihenfolge der Glieder ‚steigt'. Exponent -1: … ‚fällt'. 4) Seine Glieder sind* $h_1, h_2, h_3, …, h_s, … h_m$*. 5) s ist Position des Gliedes innerhalb des Abschnitts bei ‚steigender' Reihenfolge der Träger. 6) Bei ‚fallender' Reihenfolge gilt für Position i und Träger s:* $s + i = m + 1$.

**SATZ 3.8.** *Der geometrische Ort der Lösungsmenge der 2-dimensionalen GOLDBACH-Vermutung (2GH) ist eine Teilmenge des linearen Systems von GOLDBACH Lösungs-Rhomben:*

*1. Quadrant:* $x+y = 2 \cdot n \in (3 \cdot 2; 14) \equiv +2 \pmod{3 \cdot 2}$*. Die (x, y) werden aufgespannt von* $(H_+(m)^{+1}, H_+(m)^{-1})$*, mit* $2 \cdot n = 3 \cdot 2 \cdot (m-1)+14$*. Sie bilden eine Rhombus-Seite. Die Komponenten* $x = 3 \cdot 2 \cdot i+1$ *und* $y = 3 \cdot 2 \cdot j+1$, $i=1, 2, 3, …, i+j = m+1$, $2 \cdot n = 3 \cdot 2 \cdot (m-1)+14$*, sind symmetrisch.*

*2. Quadrant:* $x+y = 2 \cdot n \in (3 \cdot 2; 12) \equiv 0 \pmod{3 \cdot 2}$*. Die (x, y) werden aufgespannt von* $(|H_-|(m)^{+1}, H_+(H)^{-1})$*, mit* $2 \cdot n = 12+3 \cdot 2 \cdot (m-1)$*. Sie bilden eine Rhombus-Seite. Die Komponenten* $x = 3 \cdot 2 \cdot i+1$ *und* $y = 3 \cdot 2 \cdot j+1$*, mit* $i=1, 2, 3, …, i+j = m+1$, $2 \cdot n = 12 + 3 \cdot 2 \cdot (m-1)$*, sind symmetrisch.*

*3. Quadrant:* $x+y = 2 \cdot n \in (3 \cdot 2; 10) \equiv -2 \pmod{3 \cdot 2}$*. Die (x, y) werden aufgespannt von* $(H_-(m)^{+1}, H_-(m)^{-1})$*, mit* $2 \cdot n = 10+3 \cdot 2 \cdot (m-1)$*. Sie bilden eine Rhombus-Seite. Die Komponenten* $x = 3 \cdot 2 \cdot i+1$ *und* $y = 3 \cdot 2 \cdot j+1$*, mit* $i=1, 2, 3, …, i+j = m+1$, $2 \cdot n = 10 + 3 \cdot 2 \cdot (m-1)$*, sind symmetrisch.*

*4. Quadrant: Wie im 2. Quadranten, aber mit vertauschten Komponenten-Abschnitten.*

**Zu § 4. Geometrischer Ort der Lösungsmenge der Primzahlzwillinge. Verallgemeinerung**.

**Satz 4.1.** *Im Primzahlgitter ist der geometrische Ort der Lösungsmenge der Primzahlzwillinge (PT) die Winkelhalbierende des 2. (bzw. 4.) Quadranten: die Gerade* $p = (-1) \cdot q + 2$.

**Satz 4.2.** *Der geometrische Ort der Lösungen von (gPT) im Primzahl-Gitter: Die Parallelenscharen* $p = (-1) \cdot q + 3 \cdot 2 \cdot i + 2$*, i fest* $\in N$ *bzw.* $p = q + 3 \cdot 2 \cdot i$*, i fest* $\in N$*, mit: 1) verlaufen parallel zur Winkelhalbierenden des 2. (bzw. 4.) Quadranten, 2) bedecken schließlich alle vier Quadranten, 3) schneiden (auf der y-Achse) die Glieder von H.*

**Zu § 5. Geometrischer Ort der Lösungsmenge der PRACHAR Primzahl-zwillinge**. **Verallgemeinerung**.

**Satz 5.1.** *Der geometrische Ort der Lösungsmenge der PRACHAR-Primzahlzwillinge (PPT) ist die Gerade* $y = ½ \cdot (x + 1)$*, d. h. durch den Punkt (+1, +1) mit der Steigung ½. - Der geometrische Ort der Lösungsmenge der (gPPT) ist das Geraden-Bündel mit Schnitt-Punkt (+1, +1) und den Gitter-Punkten als Zweit-Punkten. - Der geometrische Ort der Lösungsmenge der total verallgemeinerten (gPPT) sind die vergleichbaren Geraden-Bündel in jedem Gitter-Punkt.*

**Zu § 6. Geometrischer Ort der Lösungsmenge der 3-gliedrigen GOLDBACH-Vermutung**.

(3GH): Lösungen der 3-gliedrigen GOLDBACHschen Vermutung.

**Satz 6.1.** *Der geometrische Ort der Lösungsmenge der (3GH) ist ein lineares System von Oktaedern. Die Oktaeder sind parallel, äquidistant mit Abstand 3·2 auf den Achsen, konzentrisch zu (1, 1, 1), die Ecken auf den Achsen. Parameter des Systems ist die Summe der Träger s; s ist fest für jeden Oktaeder. s generiert 4 ungerade Zahlen* $3 \cdot 2 \cdot s-3$, $3 \cdot 2 \cdot s - 1$, $3 \cdot 2 \cdot s + 1$ *und* $3 \cdot 2 \cdot s + 3$, $2 < s \in N$*, die vom s-ten Oktaeder dargestellt werden. Die Komponenten der Punkte sind Glieder von |H|. Die Bestimmung der Komponenten aus |H| geschieht durch den Oktaeder* $(3 \cdot 2 \cdot s_x \pm 1) + (3 \cdot 2 \cdot s_y \pm 1) + (3 \cdot 2 \cdot s_z \pm 1)$ *mit* $s = s_x + s_y + s_z$ *und* $0 \leq s_x, s_y, s_z \leq s$.

**Zu § 7. Naheliegende Vermutungen**

*Verallgemeinerung der Primzahlzwillinge. Verallgemeinerte PRACHAR Primzahlzwillinge. Verallgemeinerung der Primzahlzwillinge in den Raum. Ferner: Definition von Grundbegriffen der Geometrie allein durch Primzahlen.*

**Zu § 8. Relevanz-Intervalle**

**Satz 8.1.** *Die Primzahl-Quadrate liegen in* $H_+$.

**Satz 8.2.** $p^{2n} \in H_+$, $p^{2n+1} \in H_-$.

**Def. 8.1.** *Relevanz-Intervall, syn.* $p^2$*-Intervall: p, q* $\in P$ *und p < q und p, q aufeinander folgend in P:*

*Relevanz-Intervall* $(H_+; p^2)$ *in* $H_+ := \{ x \mid (x \in H_+) \wedge (p^2 \leq x < q^2)\}$.

*Relevanz-Intervall* $(H_-; p^2)$ *in* $H_- := \{ y \mid (y \in H_-) \wedge y = x - 2 \wedge x \in (H_+; p^2)\}$.

**Def. 8.2.** *Die relevanten Primzahlen* $p_i$ *in* $(H_+; p^2)$ *u.* $(H_-; p^2)$ *sind die Primzahlen* $p_i$ *mit* $3 < p_i \leq p$.

**Satz 8.3.** *Es ist hinreichend, in $(H_+; p^2)$ u. $(H_-; p^2)$ Vielfache als Vielfache der dortigen relevanten $p_i$ zu behandeln.*
**Def. 8.3.** *Abstand f, Länge L: Abstand f := q - p. f ist gerade. L bezeichnet die Länge (Anzahl Glieder) des Intervalls $(H_+; p^2)$. L := $(q^2 - p^2)/3\cdot 2$. -Seien $p_1 < p_2$ Primzahlzwillinge. Dann gilt: 1.) $p_1 \in H_-$ u. $p_2 \in H_+$, 2.) Das Glied zwischen $p_1$ u. $p_2$ wird von 3 geteilt, 3.) Nur bei Primzahlzwillingen ist L < p; sonst ist p < L.*
**Satz 8.4.** *p aufeinander folgende Glieder von Z, N, H, $H_+$, $H_-$ lassen ein vollständiges Restsystem modulo p.*
**Satz 8.6.** *1.) Die Vereinigung der Relevanz-Intervalle ist IHI - abgesehen vom Glied -1, 2 und 3. 2.) $(H_+; p^2)$ u. $(H_-; p^2)$ besitzen die gleiche Länge L 3.) $(H_+; p^2)$ enthält mindestens ein Vielfaches, u. zwar $p^2$. 4.) $(H_+; p^2)$ u. $(H_-; p^2)$ haben dieselben relevanten Primfaktoren. 5.) Die Primfaktor-Zerlegung eines Vielfachen im Rel.-Int. enthält mindestens einen relevanten Primfaktor. 6.) Die Vielfachen im Rel.-Int. sind Vielfache von welchen der relevanten Primfaktoren. 7.) In jedem weiteren Rel.-Int. tritt genau ein relevanter Primfaktor hinzu, u. zwar p.*
**Satz 8.7.** *p generiert die arithmetische Folge $(\pm 3\cdot 2\cdot p; p^2) \subset H$ der Vielfachen von p.*
**Satz 8.8.** *Alle Glieder von $(H_+; p_i^2)$ oder $(H_-; p_i^2)$, die in einer arithmetischen Folge $(\pm 3\cdot 2\cdot p_i; p_i^2)$ vorkommen, sind Vielfache. Alle anderen sind Primzahlen.*
**Satz 8.9.** *Die Vielfachen von p, sofern sie zu $(\pm 3\cdot 2\cdot p; p^2)$ gehören, erscheinen abwechselnd in den Ästen von H.*
*$H_+$ enthält $-5\cdot H_-$, $+7\cdot H_+$, $-11\cdot H_-$, $+13\cdot H_+$, $-17\cdot H_-$, ... $H_-$ enthält $-5\cdot H_+$, $+7\cdot H_-$, $-11\cdot H_+$, $+13\cdot H_-$, $-17\cdot H_+$, ...*

**Zu § 9. Primzahlen in Relevanz-Intervallen.**
**Satz 9.1.** *Es gibt unendl. viele Relevanz-Interv. der Gestalt $(H_+; p^2)$ und es gibt unendlich viele Relevanz-Intervalle der Gestalt $(H_-; p^2)$.*
**Satz 9.2.** *Es gibt unendlich viele Relevanz-Intervalle in $H_+$, die Primzahlen enthlaten, und es gibt unendlich viele Relevanz-Intervalle.*
**Satz 9.3.** *Es gibt in $(H_+; p^2)$ etwa gleich viele Vielfache wie in $(H_-; p^2)$.*
**Satz 9.4.** *Wenn es in $(H_+; p^2)$ Primzahlen gibt, dann auch in $(H_-; p^2)$ - und umgekehrt - und zwar etwa gleich viele.*
**Satz 9.5.** *Es gibt unendlich viele ‚parallele' Relevanz-Intervalle, die Primzahlen enthalten.*

**Zu § 10. Unendlichkeit der Menge der Primzahlzwillinge.**
**Satz 10.1.** *Wenn $\#v_-v_+$** > 0 ist, so gilt auch $\#p_-p_+ > 0$ - und umgekehrt.*
**Satz 10.2.** *Es gibt ∞ viele Primzahlzwillinge.*
**Satz 10.3.** *Ist in parallelen Relevanz-Intervallen #v > #p, so gibt es dort vv- und vp- und pv- und pp-Paare.*
**Satz 10.4.** *Schritte 1 bis 20 funktionieren bei jeder Primzahl $p \geq 23$. Dadurch ist eine ‚Struktur' der Primzahlverteilung beschrieben.*

**Zu § 11. Bestätigung der GOLDBACH-Vermutung**
**Definition 11. 1.** *Ein Achsen-Abschnitt $H_+(A)^{+1}$ ist ein Anfangsstück $\{h_1, h_2, h_3, \ldots, h_A\}$ von $H_+$ der Länge A, das bei Exponent +1 steigend durchlaufen wird, bei Exponent -1 fallend. Analog in den drei andereni Quadranten.*
**Definition 11.2.** *Eine Rhombusseite $(H_+(A)^{+1}, H_+(A)^{-1})$ wird aufgespannt von den Achsen-Abschnitten $H_+(A)^{+1}$und $H_+(A)^{-1}$). Ihre Punkte sind die Summen der aufspannenden Koordinaten. Analog bei den drei anderen Seiten des Rhombus.*
**Satz 11.1.** *#pp = [A -(#pv + #vp) - (#v - #p)] / 2. #vv = [A -(#pv + #vp) - (#v - #p)] / 2 + [#v - #p]. A ist die Länge des Komponenten-Abschnitts.*
**Satz 11.2. Bestätigung der GOLDBACH-Vermutung.**
*Die Lösungen der* GOLDBACH-Vermutung *(2GH) sind darstellbar als Punkte im Primzahl-Gitter. Der geometrische Ort der Lösungsmenge der (2GH) im Primzahl-Gitter ist ein System von achsen-symmetrischen Rhomben. Es gilt: 1.) Jede gerade Zahl besitzt mindestens eine, im Allgemeinen mehrere (2GH)-Lösungen. 2.) Alle (2GH)-Lösungen, die eine bestimmte gerade Zahl darstellen und nur diese, befinden sich auf einer ihnen umkehrbar eindeutig zugeordneten Rhombus-Seite. 3.) Die Gestalt der (2GH)-Lösungen ist durch ihre Rhombus-Seite bestimmt. 4.) Man kann die Lösungen anhand der Rhombus-Seite explizit angeben. 5.) Es gibt für die zu einer geraden Zahl gehörenden Lösungen - außer ihrer Rhombus-Seite - keine geschlossene Darstellung. – Die (2GH) ist zutreffend.*

## Literatur

[1] PRACHAR, Karl. Primzahlverteilung. Reprint. Berlin, Heidelberg, New York: Springer, 1957.
[2] SCHOLZ, Arnold, u. SCHOENEBERG, Bruno. Einführung in die Zahlentheorie. Sammlung Göschen, Band 1131. Berlin: Walter de Gruyter & Co., 1955.
[3] SCHWARZ, Wolfgang. Einführung in Methoden u. Ergebnisse der Primzahltheorie. B·I·Hochschul-Taschenbücher 278/278a*. Mannheim, Wien, Zürich: Bibliograph. Institut, Hochschultaschenbuch-Verlag, 1969.
[4] REINHARDT, Fritz, u. SOEDER, Heinrich. dtv Atlas zur Mathematik, Tafeln und Texte. Band I, Grundlagen, Algebra und Geometrie. Deutscher Taschenbuch Verlag. 9. Auflage, 1991.
[5] HUBBARD, R. L. The Factor Book. Hilton Management Services Ltd., Lytham St. Annes, Lancashire, England, 1975.